\documentclass[twocolumn,final]{svjour3}           % twocolumn
\smartqed  % flush right qed marks, e.g. at end of proof
\usepackage{upgreek}
\usepackage{amsmath,amssymb,graphicx,color}
\usepackage{epstopdf}
\usepackage{algorithm,algorithmic}
\usepackage{cases}
\usepackage{mathtools, cuted}
\usepackage{siunitx}
\usepackage{hyperref}
\usepackage{mathrsfs}

\usepackage[varg]{txfonts}
\newcommand{\phii}{\phiup}                     %% porosité
\newcommand{\R}{\mathbb{R}}                  %% ensemble reel
\newcommand{\G}{\textbf{g}}                  %% terme de gravite
\newcommand{\N}{\mathbb{N}}                  %% ensemble entier
\newcommand{\bta}{\beta}                     %% composant
\newcommand{\Ka}{{\bf K}}                    %% permeabilite absolue

\newcommand{\J}{\vec{J}}                     %% flux
\newcommand{\V}{\textbf{V}}                  %% vitesse de Darcy
\newcommand{\Dlh}{D_l^h}                     %% terme de diffusion
\newcommand{\ctea}{\mathcal{C}}              %% constante a
\newcommand{\X}{\mathscr{X}}                 %% fraction massique
\newcommand{\XX}{\mathcal{X}}                 %% fraction massique
\newcommand{\sa}{S} %% saturation S

\newcommand{\dt}{\partial_{t}}
\newcommand{\di}{\nabla\cdot}
\newcommand{\dd}{\text{d}}
\newcommand{\grad}{\nabla}

\newcommand{\ip}[2]{\langle{#1}, {#2}\rangle}
\newcommand{\mat}[1]{\mathbf{{#1}}}

\newcommand{\Npc}{{N}_{\text{PC}}}
\newcommand{\Mpc}{{M}_{\text{PC}}}
\newcommand{\Ns}{{N}_{\text{s}}}
\newcommand{\Cpar}{\mathcal{C}_{\text{par}}}
\newcommand{\Cqoi}{\mathcal{C}_{\text{qoi}}}
\newcommand{\Nqoi}{{{N}_{\text{qoi}}}}
\newcommand{\Mqoi}{{M}_{\text{qoi}}}
\newcommand{\Np}{{{N}_{\text{p}}}}
\newcommand{\np}{{{n}_{\text{p}}}}

\newcommand{\Nord}{{N}_{\text{ord}}}
\newcommand{\PC}{\text{PC}}

\renewcommand{\vec}[1]{{\mathchoice
                     {\mbox{\boldmath$\displaystyle{#1}$}}
                     {\mbox{\boldmath$\textstyle{#1}$}}
                     {\mbox{\boldmath$\scriptstyle{#1}$}}
                     {\mbox{\boldmath$\scriptscriptstyle{#1}$}}}}

\DeclareMathOperator*{\argmin}{arg\,min}

\newcommand{\dm}{\ell}
\newcommand{\Dalfabta}{D_\alpha^\beta}    %%  diffusion term

\newcommand{\erel}{\text{e}_{\text{rel}}}

\newcommand{\sj}{\mathfrak{s}_j}

\newcommand{\normLomega}[1]{\left\|#1\right\|_{L^2(\paramspace)}}
\newcommand{\norm}[1]{\left\|#1\right\|}
\newcommand{\fn}{f_{{{\Nqoi}\!}}}
\newcommand{\gn}{g_{{{\Mqoi}\!}}}
\newcommand{\Sn}{S_{{{\Nqoi}\!}}}
\newcommand{\Qn}{Q_{{{\Nqoi}\!}}}

\newcommand{\fPC}{\fn^{{\text{PC}}}}
\newcommand{\gPC}{\gn^{{\text{PC}}}}
\newcommand{\SPC}{\Sn^{\PC}}
\newcommand{\QPC}{\Qn^{\PC}}

\newcommand{\coma}{\eta}
\newcommand{\comb}{\tilde{\coma}}
\newcommand{\Cov}[1]{\text{Cov}\{#1\}}
\newcommand{\paramspace}{\Theta}
\newcommand{\mymu}{\!\upmu}

\newcommand*{\transymb}{{\mkern-1.5mu\mathsf{T}}}
\renewcommand{\top}{\transymb}
\newcommand*{\Scale}[2][4]{\scalebox{#1}{\ensuremath{#2}}}%

\newcommand{\bigcomma}{\Scale[1.5]{,}}

\graphicspath{{figures/},{tikzpictures/}} %Setting the graphicspath
\allowdisplaybreaks

\begin{document}

\title{
Structure exploiting methods for fast uncertainty quantification in multiphase flow 
through heterogeneous media
}

%\subtitle{Do you have a subtitle?\\ If so, write it here}

%\titlerunning{Short form of title}        % if too long for running head

\author{
Helen Cleaves             \and
Alen Alexanderian         \and
Bilal Saad
}

\institute{H.~Cleaves \at
              Department of Mathematics,
              North Carolina State University, \\
              Box 8205,
              Raleigh, NC 27695. \\
              \email{hlcleave@ncsu.edu}           
           \and
           A.~Alexanderian \at
              Department of Mathematics,
              North Carolina State University, \\
              Box 8205,
              Raleigh, NC 27695.\\
              \email{alexanderian@ncsu.edu}           
           \and
           B.~Saad \at
              Ecole Centrale de Nantes,\\
              1, rue de la No\'e, 44321 Nantes, France. \\
              \email{bisaad@yahoo.fr}
}

\date{Received: date / Accepted: date}
% The correct dates will be entered by the editor

\maketitle

\begin{abstract}
We present a computational framework for dimension reduction and surrogate
modeling to accelerate uncertainty quantification in 
computationally intensive models with high-dimensional inputs and
function-valued outputs. Our driving application is multiphase flow in saturated-unsaturated
porous media in the context of radioactive waste storage. For fast input dimension
reduction,  we utilize an approximate global sensitivity measure, for
function-valued outputs, motivated by ideas from the active subspace methods.
The proposed approach does not require expensive gradient computations.  We
generate an efficient surrogate model by combining a truncated
Karhunen-Lo\'{e}ve (KL) expansion of the output with polynomial chaos
expansions, for the output KL modes, constructed in the reduced parameter
space.  We demonstrate the effectiveness of the proposed surrogate modeling
approach with a comprehensive set of numerical experiments, where we consider a number
of function-valued (temporally or spatially distributed) QoIs.

\keywords{Uncertainty quantification \and surrogate models \and dimension reduction \and
multiphase flow \and sensitivity analysis \and spectral representations}

\subclass{65C20 \and 65C50 \and 65D15 \and 76S05 \and 35Q86}

\end{abstract}

% This is all for the toc
%\setcounter{secnumdepth}{3}
%\setcounter{tocdepth}{3}
%\tableofcontents

\section{Introduction}
\label{intro}
Low permeability argillites are considered as suitable host rocks for
underground radioactive waste storage to retain radionuclides locally. However,
hydrogen gas produced by corrosion of steel engineered barriers can represent a
threat to the installation safety. A significant impact of this production is
the overpressurization of hydrogen around alveolus leading to opening fractures
in the surrounding host rock and inducing groundwater flow and transport of
radionuclides outside of the geological repositories. This problem renews the
mathematical interest in the equations describing multiphase multicomponent
flows through porous media, within the present context.  An important aspect of
improving the prediction fidelity of such models is to account for the various
sources of uncertainty in the governing equations.

Performing uncertainty analysis on the models under study using a direct Monte
Carlo sampling approach is infeasible. This is due to the high cost of model
simulations and the need for a large number of such simulations.  Therefore,
there is a need for quick-to-evaluate surrogate models that accurately capture
the underlying physics and statistical properties of the quantities of
interest~(QoIs). Surrogate modeling, however, is a formidable task for the applications
considered in the present work. Models describing flow through porous media
exhibit distinct challenges with regards to uncertainty quantification and surrogate modeling
including expensive simulations, high-dimensional uncertain parameters, and
function-valued outputs.  Addressing these challenges effectively requires
understanding and exploiting the problem structure. To this end, we propose a
framework that deploys a sensitivity analysis approach to reduce the
dimensionality of the input parameter and utilizes the spectral properties of 
the output QoI to generate an efficient surrogate model.

\begin{figure*}[h]\centering
\includegraphics[width=.9\textwidth]{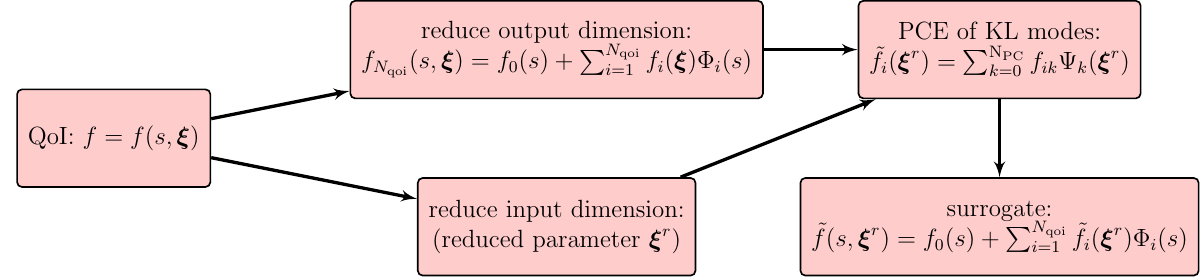}
  \caption{A schematic of the proposed \emph{bispectral} surrogate modeling approach.}
\label{fig:diag}
\end{figure*}

\textbf{Related work.}
The modeling of underground radioactive waste storage involves simulating
the coupled transport of multiphase multicomponent flow in porous medium.
Equations governing this type of flow in porous media are nonlinear and
involve simulation of complex phenomena such as the appearance and the
disappearance of the gas phase leading to the degeneracy of the equations
satisfied by the saturation. There have been significant research efforts dealing with
mathematical and numerical models for simulating the transport migration of
radionuclides. The articles~\cite{Bourgeat2012,Bourgeat2009} present test-cases
and set up benchmark examples to address some of the specific problems encountered
when numerically simulating gas migration in underground nuclear waste repositories.
In~\cite{Angelini2011,Bourgeat2009,Neumann2012} different choices of primary
variables have been proposed to tackle the degeneracy of the equations satisfied
by the saturation. In~\cite{Saad2020}, the authors study a compressible
and partially miscible phase flow model in porous media, applied to gas migration in
an underground nuclear waste repository in the case where the velocity of the
mass exchange between dissolved hydrogen and hydrogen in the gas phase is
supposed finite. Also presented is a numerical scheme based on a two-step
convection/diffusion-relaxation strategy to simulate the non-equilibrium model.
There have also been efforts to quantify uncertainty in models of multiphase
flow~\cite{ChristieDemyanovErbas06,Najm09,NamhataOladyshkinDilmore16,SaadAlexanderianPrudhommeEtAl18,SeverinoLevequeToraldo19,XiaoOladyshkinNowak19}.

The tools from uncertainty quantification that are relevant to the
present work include global sensitivity analysis (GSA) and surrogate modeling.
GSA provides insight into how uncertainties in
model parameters influence model outputs by identifying the input parameters a
QoI is sensitive to. This increases overall understanding of the
underlying physics and guides parameter dimension reduction.  The Sobol'
indices~\cite{Sobol01}, derivative-based global sensitivity measures~(DGSMs)
\cite{KucherenkoRodriguez09,KucherenkoIooss17,SobolKucherenko09}, and active
subspace methods~\cite{Constantine15,ConstantineDiaz17} are examples of GSA tools
widely used in practice.  These concepts were originally conceived for
scalar QoIs. Recent works such
as~\cite{AlexanderianGremaudSmith20,CleavesAlexanderianGuyEtAl19,ZahmConstantine2020}
generalize standard GSA tools to the case of vector- and function-valued QoIs.
In particular,~\cite{GamboaJanonKlein14,AlexanderianGremaudSmith20} concern
variance-based GSA using Sobol' indices for such QoIs.  The
article~\cite{CleavesAlexanderianGuyEtAl19} studies DGSMs for function-valued QoIs.
A generalization of active subspace methods for vectorial outputs is presented
in \cite{ZahmConstantine2020}.

For expensive-to-compute QoIs calculating GSA measures such as Sobol' indices is
computationally expensive. A common method for mitigating the computational
cost is to construct a cheap-to-evaluate surrogate model for the QoI and then apply
GSA techniques to the surrogate. For example, polynomial
chaos expansions~(PCEs) have been a popular approach for accelerating the
computation of Sobol' indices; see,
e.g.,~\cite{AlexanderianWinokurSrajEtAl12,BlatmanetAl2010,Crestaux09,Sudret08}.
Surrogate model construction, however, is itself a computationally challenging
task, especially in the case of models with high-dimensional input parameters.  For
such models it is also possible to use a multilevel approach: initial parameter
screening can be performed using cheap, but less precise, tools and further
dimension reduction is performed through more rigorous methods such as a
variance-based analysis using accurate surrogate models constructed in a
reduced-dimensional parameter space; see e.g.,~\cite{HartGremaudDavid19}.

For function-valued QoIs, a straightforward approach is to compute surrogate
models for every grid point in a discretized computational domain.  This
approach, however, can be inefficient and ignores an important problem
structure---the low-rank structure of the output.  Specifically, in many
applications, function-valued QoIs can be represented via a spectral
representation, such as a Karhunen--Lo\'{e}ve expansion~(KLE), with a small
number of terms.  This problem structure can be exploited for surrogate
modeling: instead of approximating a field quantity at every point in a
computational grid, one can approximate a few dominant modes of the output QoI.
Such surrogate  models can also be used to accelerate GSA methods; see
e.g.,~\cite{AlexanderianGremaudSmith20,CleavesAlexanderianGuyEtAl19,GuyetAl19,LiIskandaraniLeHenaff16}.

\textbf{Our approach and contributions.}
%How we explaining our contribution, he said will be simpler that we said we
%generalize the SRC but we are not. They are related but not the same things...
%Should we write remarks: Normalized function activity score correspond to the
%function Sobol indices for the linearized model in the case where the standard
%deviation is one...

In the present work, we seek to construct surrogate models for fast analysis of
computationally intensive models with high-dimensional parameters and
function-valued QoIs. We consider QoIs of the form
\begin{equation*}
f = f(s, \vec\xi), \quad s \in \X, \, \vec\xi \in \Theta,
\end{equation*}
where $\Theta \subseteq \R^{\Np}$ is the uncertain parameter domain and $\X$ is
compact subset of $\R^d$, with $d \in \{1, 2, 3\}$. In practice, $s$ can
represent a spatial or temporal point.  Our focus in the present work is models
of flow in porous media, and $f(s, \vec\xi)$ is an observable in a multiphase
flow problem.  Our approach identifies and exploits low-dimensional structures
in both input and output spaces. Specifically, we rely on approximate GSA
measures for fast input parameter screening and utilize low-rank spectral
representations of output fields.

We propose a fast-to-compute screening metric that utilizes ideas from active
subspaces~\cite{Constantine15} and derivative-based GSA for functional 
outputs~\cite{CleavesAlexanderianGuyEtAl19}
to perform parameter dimension reduction.
%
%Specifically, we use a global linear regression model, for
%approximating activity scores of function-valued QoIs. It is worth noting
%that linear regression analysis has been widely 
%used for obtaining sensitivity measures, such as standardized regression 
%coefficients, 
%for scalar QoIs; see e.g.,~\cite{Helton93,ConstantineDiaz17}. 
%Our work extends such analyses to the case of function-valued QoIs
%in a systematic manner.
%
The proposed screening metrics do not require gradient computation
in the parameter space. This makes the proposed methods applicable to a broad
class of problems involving complex physics systems for which adjoint solvers,
which are essential for gradient computation in high dimensions, are not
necessarily available.  

Following parameter screening, we combine two different
spectral approaches---KLEs and PCEs---to generate an efficient surrogate model
in a reduced-dimensional uncertain parameter space.  The overall surrogate
model constructed takes the form,
\begin{equation*}
\label{equ:kle_surrogate}
\fPC(s, \vec\xi^r) = f_0(s) +
  \sum_{i=1}^\Nqoi f_i^{\PC}(\vec\xi^r) \Phi_i(s),
\end{equation*}
where $\Phi_i$'s are orthogonal basis functions in $L^2(\X)$ obtained from a KLE of $f(s,
\vec\xi)$ and $f_i^{PC}$ are approximate KL modes as functions of a
reduced-dimensional parameter vector $\vec\xi^r\subseteq\R^\np$; these KL modes are
represented by PCEs,
\[
  f_i^{\PC}(\vec\xi^r) = \sum_{k=0}^{\Npc} f_{ik}
\Psi_k(\vec\xi^r),
\]
where $\Psi_k$'s are a basis consisting of multivariate orthogonal polynomials
in $L^2(\Theta)$ and $\Npc$ is specified based on the choice of truncation strategy. 
Thus, the overall surrogate model can be expressed as
\begin{equation}\label{equ:overall_surrogate}
f(s, \vec\xi) \approx f_0(s) + \sum_{i=1}^\Nqoi \sum_{k=0}^{\Npc} f_{ik} \Psi_k(\vec\xi^r) \Phi_i(s).
\end{equation}
We refer to the class of surrogate models of the
form~\eqref{equ:overall_surrogate} as \emph{bispectral surrogates} due to the
use of spectral representations in $L^2(\X)$ and $L^2(\Theta)$. In
Figure~\ref{fig:diag}, we provide a schematic of the proposed bispectral
surrogate modeling framework. We point out that the proposed
approach is non-intrusive and requires only the ability to
evaluate the governing model at a sample of uncertain inputs.
See Section~\ref{sec:methods} for details.

While computing a surrogate model from a truncated KLE by replacing the KL
modes with PCEs (or other surrogates) is not new, see
e.g.,~\cite{AlexanderianGremaudSmith20,GuyetAl19,LiIskandaraniLeHenaff16}, we
build upon this approach by including a gradient-free input dimension reduction
approach as a first step.  This enables the PCEs for the KL modes to be built
in a lower-dimensional space.  Thus, a major contribution of this article is a
synergy of known techniques combined with a novel input dimension reduction
strategy to furnish an integrated surrogate modeling approach.  We also provide
a detailed computational procedure for the proposed framework, making the
present work a self-contained guide.  We elaborate our approach on an intricate
multiphase multicomponent flow model for which a comprehensive presentation is
also given. In our numerical results, we implement the proposed approach for
both spatially- and temporally-varying QoIs.  Additionally, a variety of
statistical studies are conducted with the constructed bispectral surrogate.
These tests are intended to showcase the versatility of the surrogate model and
explore the physical phenomenon under study.  In particular, we perform model
predictions, compute variance-based global sensitivity indices, and study
statistical model response behavior. In addition to demonstrating the
effectiveness of the proposed strategy,  our computational results provide
valuable insight regarding the response of complex porous media flow models to
uncertainties in material properties.

\textbf{Article overview.}
In section~\ref{sec:model} we present a
detailed overview of the multiphase multicomponent flow model that is central
to the present work. We also provide a description of our choice of numerical
solver for the governing equations.  In section~\ref{sec:model_uncertainty}, we
discuss modeling the uncertainties in material properties, as well as give a
brief explanation of the model response and relevant QoIs. We supply a concise
overview of KLEs, PCEs, and bispectral surrogates in section~\ref{sec:spectral}. In
section~\ref{sec:methods} we provide a detailed framework, including
algorithms, for the proposed dimension reduction and surrogate modeling
approach. Our computational results are presented in
section~\ref{sec:numerics}. Finally, we provide closing comments in
section~\ref{sec:conc}.

\section{Model Description}
\label{sec:model}
\subsection{Mathematical formulation of the continuous problem}

Here we state the physical model used in this work. We consider a porous medium
saturated with a fluid composed of two phases, liquid ($l$) and gas ($g$),
and a mixture of two components, water ($w$) and hydrogen ($h$).
The spatial domain $\Omega$ is a bounded open subset of $\R^\dm$ ($\dm = 1, 2$, or $3$)
and the problem is considered in the time interval $[0, T_f]$, where $T_f > 0$
is the final time. To define the physical model, we write the
\emph{mass conservation} of each component in each phase
\begin{align}
  \phii \partial_{t}(\rho_{l}^w \sa_{\! l} &+ \rho_{g}^w \sa_{\! g})\notag \\
  &+\di(\rho_{l}^w \V_{l} + \rho_{g}^w \V_{g} + \J_l^{w} + \J_g^{w}) =
  f^{w},
  \label{equ:sys_eqA}\\
  \phii\partial_{t}(\rho_{l}^h \sa_{\! l} &+ \rho_{g}^h \sa_{\! g}) \notag\\
  &+ \di(\rho_{l}^h \V_{l} + \rho_{g}^h \V_{g} + \J_l^{h} + \J_g^{h}) = f^{h},
  \label{equ:sys_eqB}
\end{align}
where  $\phii(x)$ is the given porosity of the medium, $\sa_{\!\alpha}(t,x)$ the
saturation of the phase $\alpha \in \{l, g\}$, with the two saturations summing to
one.
Also, $p_\alpha(t,x)$ is the pressure of the phase $\alpha$,
$\rho_\alpha^{\bta}$ is the density of the component $\bta \in \{w, h\}$
in the phase $\alpha$, and $\rho_\alpha=\rho_\alpha^{h}+ \rho_\alpha^{w}$ is
the density of the phase $\alpha$. The velocity of each fluid, $\V_\alpha$
is given by Darcy's law
\begin{equation*}
  \V_{\alpha}= -{\bf K} \frac{ k_{r_\alpha}(\sa_{\! \alpha})}{\mu_{\alpha}}
  \big(\nabla p_{\alpha}-\rho_{\alpha}(p_{\alpha})\textbf{g}\big),
\end{equation*}
where $\Ka(x)$ is the intrinsic (given) permeability tensor of the porous
medium, $k_{r_\alpha}$ the relative permeability of the $\alpha$-phase,
$\mu_\alpha$ the constant $\alpha$-phase's viscosity, $p_\alpha$ the
$\alpha$-phase's pressure, and $\G$, the gravity vector.  For further details of
the model we refer to the presentation of the benchmark~\cite{Bourgeat2012,Bourgeat2009}.
Following the Fick's law, the diffusive flux of a component $\bta$ in the phase
$\alpha$ is given by
\begin{equation*}\label{eq:flux_diffusif}
  \J_\alpha^\beta = - \phii \sa_{\!\alpha} \rho_\alpha \Dalfabta \grad \XX_\alpha^\beta,
\end{equation*}
where coefficient $\Dalfabta$ is the Darcy scale molecular diffusion
coefficients of $\bta$-component in $\alpha$-phase and $\XX_\alpha^\beta
= {\rho_\alpha^\beta}/{\rho_\alpha}$ is the component $\beta$ molar
fraction in phase $\alpha$. Diffusive fluxes satisfy
$\sum_{\beta}^{}\J_\alpha^\beta = 0$ for each $\alpha$.

The capillary pressure law, which links the jump of pressure of the two phases
to the saturation, is
\begin{equation*}
  p_c(\sa_{\!l}) = p_g - p_l.
\end{equation*}
This function is decreasing ($\frac{\dd p_c}{\dd \sa_{\!l}}(\sa_{\!l})<0 \text{ for all
} \sa_{\!l}\in[0,1]$), and satisfies $p_c(1)=0$.

In the present work, the water is supposed only present in the liquid phase (no
vapor of water due to evaporation). Thus, \eqref{equ:sys_eqA}--\eqref{equ:sys_eqB}
could be rewritten as
\begin{align}
  \phii\dt(\sa_{\!l}\,\rho_l^{w} ) &+ \di(\rho_l^{w} \V_l)\notag\\
  &+ \di(\phii \sa_{\!l}\, \rho_l\Dlh\nabla X_l^h) = f^w, \label{eq:water} \\
  \phii\dt(\sa_{\!l}\,\rho_l^{h} + \sa_{\!g}\,\rho_g^{h} )
   &+ \di(\rho_l^{h} \V_l + \rho_g^{h} \V_g)
\notag
\\
  &-\di(\phii \sa_{\!l}\, \rho_l\Dlh\nabla X_l^h) = f^h.
\label{eq:hydrogen}
\end{align}
The system \eqref{eq:water}--\eqref{eq:hydrogen} is not complete; to close the
system, we use the ideal gas law and the Henry's law
\begin{equation}\label{law}
  \rho_g^h = \frac{M^h}{R T}p_g\quad \text{and} \quad \rho_l^h = M^h H^h p_g,
\end{equation}
where the quantities $M^h$, $H^h$, $R$ and $T$ represent respectively the molar
mass of hydrogen, the Henry's constant for hydrogen, the universal constant of
perfect gases and $T$ the temperature.
By these formulation, the system \eqref{eq:water}--\eqref{eq:hydrogen} is
closed and we choose the liquid pressure and the density of dissolved hydrogen
as unknowns.  From \eqref{law}, the Henry's law combined to the ideal gas law,
to obtain that the density of hydrogen gas is proportional to the density of
hydrogen dissolved
\begin{equation*}\label{eq:relation_densite}
  \rho_g^h = \ctea \rho_l^h \quad\text{where}\quad \ctea = \frac{1}{H_h R T} = 52.51.
\end{equation*}

Note that the density of water $\rho_l^w$ in the liquid phase is constant and
from the Henry's law, we can write
\begin{equation*}
  \rho_l\nabla \XX_l^h = \XX_l^w\nabla \rho_l^h.
\end{equation*}
Then the system \eqref{eq:water}--\eqref{eq:hydrogen} can be written as
\begin{align}
  \phii \partial_{t}\left(\sa_{\!l}\,\rho_l^w \right) &+
  \di\left(\rho_l^w \V_{l}\right)\notag \\
  &+ \di\left(\phii \sa_{\!l}\, \XX_l^w \Dlh \nabla \rho_l^h \right)
  = f^w,  \notag\\ %\label{eq:water1} \\
 \phii \partial_{t}\left(m(\sa_{\!l}) \rho_l^h \right) &+
  \di\left(\rho_l^h\V_{l} + \ctea \rho_l^h\V_{g}\right) \notag\\
  &- \di\left(\phii \sa_{\!l}\, \XX_l^w \Dlh \nabla \rho_l^h \right)=f^h, \notag\\
%\label{eq:hydrogen1}
  \end{align}
where $m(\sa_{\! l})=\sa_{\! l} + \ctea \sa_{\! g}$.\\

A van Genuchten-Mualem model with the parameters $n$, $\sa_{\!\alpha r}$ and $p_r$
as given in Table \ref{tab:nomvals}~(left) is used for the relative permeabilities and
capillary pressure:
\[
\begin{aligned}
  p_c(\sa_{\! le})     &= p_r \left( \sa_{\! le}^{-1/\upsilon} - 1  \right)^{1/n} ,\\
  k_{r_l}(\sa_{\! le}) &= \sqrt{\sa_{\! le}} \left( 1 - \left( 1 - \sa_{\! le}^{1/\upsilon} \right)^\upsilon   \right)^2, \\
  k_{r_g}(\sa_{\! le}) &= \sqrt{1 - \sa_{\! le}}  \left( 1 - \sa_{\! le}^{1/\upsilon} \right)^{2\upsilon},
\end{aligned}
\]
with the effective saturation
\begin{equation*}
\sa_{\! le} = (\sa_{\!l} - \sa_{\! lr})/(1 - \sa_{\! lr} - \sa_{\! gr}),
\end{equation*}
where  $\sa_{\! lr}$ and $\sa_{\! gr}$ are the liquid and gas
residual saturations, respectively, and $\upsilon = 1 - 1/n$.

%
% Numerical method
%
\subsection{Numerical solver}
As is well known, the modeling of underground radioactive waste storage
involves simulation of complex phenomena such as the appearance and the
disappearance of the gas phase leading to the degeneracy of the equations
satisfied by the saturation. This is mainly due to the migration of gas
produced by the corrosion of nuclear waste packages within a
complex heterogeneous domain. To overcome this difficulty, an important
consideration, in the modelling of multiphase flow with mass exchange between
phases, is the choice of the primary variables that define the thermodynamic
state of the system. Different choices of primary variables have been proposed
~\cite{Angelini2011,Bourgeat2009,Neumann2012}.  In this article, we consider
pressure of the liquid phase and density of dissolved hydrogen the primary
unknowns in the multiphase flow system.
A cell-centered finite volume scheme is used for the space discretization
and an implicit Euler scheme for the temporal discretization. The nonlinear
system is solved with a fixed point method.

In this section, we present a numerical study dedicated to understanding the
computational issues caused by gas phase appearance produced by injecting of
hydrogen in a one-dimensional homogeneous porous domain. 
We consider a domain that is fully saturated with water. 
This numerical study is inspired by the MoMaS benchmark on multiphase flow in
porous media~\cite{Bourgeat2012}.

\subsection{Numerical experiment}
We consider a one-dimensional domain with the benchmark setup described
in~\cite{Bourgeat2012}. The spatial domain $\Omega$ is the interval $(0, L)$,
with $L = 200$ meters, and the final simulation time is $T_f = 10^6$ years.
The parameters for porous medium, fluid characteristics, and initial and
boundary conditions are presented in~\cite{Bourgeat2012} and summarized in
Table~\ref{tab:nomvals}.
\begin{table*}\centering
\caption{Left: parameter values for the porous medium and fluid characteristics
used in test case 1. Right: parameter values for domain size, boundary and
initial conditions, total injection time and total simulation time.
}
\begin{tabular}{cc}
\begin{tabular}{llllll}
\hline
{\it Parameter}            & {\it Value}     & & {\it Parameter}                            & {\it Value}          \\
\hline
$\phii$ [-]                 & 0.15            & & $\Dlh$  [$\mathrm{m^2\cdot s^{-1}}$]       & $3 \times 10^{-9}$   \\
${\bf K}$ [$\mathrm{m^2}$] & $5 \times 10^{-20}$    & & $\mu_l$ [$\mathrm{Pa \cdot s}$]            & $1\times 10^{-3}$    \\
$p_r$ [$\mathrm{Pa}$]      &  $2\times 10^6$ & & $\mu_g$ [$\mathrm{Pa \cdot s}$]            & $9\times 10^{-6}$    \\
$n$ [-]                    & $1.54$          & & $H^h$ [$\mathrm{mol.Pa^{-1}.m^{-3} }$]     & $7.65\times 10^{-6}$ \\
$\sa_{lr}$ [-]               & $0.4$           & & $M^h$ [$\mathrm{Kg \cdot mol^{-1} }$]      & $2\times 10^{-3}$    \\
$\sa_{gr}$ [-]               & $0$             & & $\rho_l^w$ [$\mathrm{Kg \cdot mol^{-3} }$] & $ 10^{3}$            \\
\hline
\end{tabular}
&
\begin{tabular}{llllll}
\hline
{\it Parameter}                                & {\it Value}         \\
\hline
$L$ [$\mathrm{m}$]                        & 200            \\
$q_h$ [$\mathrm{kg/m^2/year}$]   & $5.57\times 10^{-6}$ \\
$p_{init}$ [$\mathrm{Pa}$]              & $10^6$   \\
$T_{inj}$ [$\mathrm{years}$]           & $5 \times 10^{5}$ \\
$T_{f}$ [$\mathrm{years}$]              & $10^{6}$ \\
\hline
\end{tabular}
\end{tabular}
\label{tab:nomvals}
\end{table*}

Initial conditions are uniform over the whole domain with pure liquid water
at fixed liquid pressure and no hydrogen present,
\begin{equation*}
p_l (0, x)= p_{init} \quad \text{and} \quad \rho_l^h (0,x) = 0,
\quad x \in \Omega.
\end{equation*}
For boundary conditions, a constant flux of hydrogen and zero water flow rate
were imposed on the left boundary
\[
\begin{aligned}
\rho_l^{w} V_l - J_l^h &= 0,\\
\rho_l^{h} V_l + \rho_g^{h} V_g + J_l^h &=
\begin{cases}
q_h \quad 0 \leq t  \leq T_\text{inj},\\
0 \quad t > T_\text{inj}.
\end{cases}
%q_h \chi_{[0,T_{\text{inj}}]}(t),
\end{aligned}
\]
On the right boundary, Dirichlet boundary conditions the same as the initial
conditions are imposed.

\begin{figure}[htp]
 \centering
  \includegraphics[width=.4\textwidth]{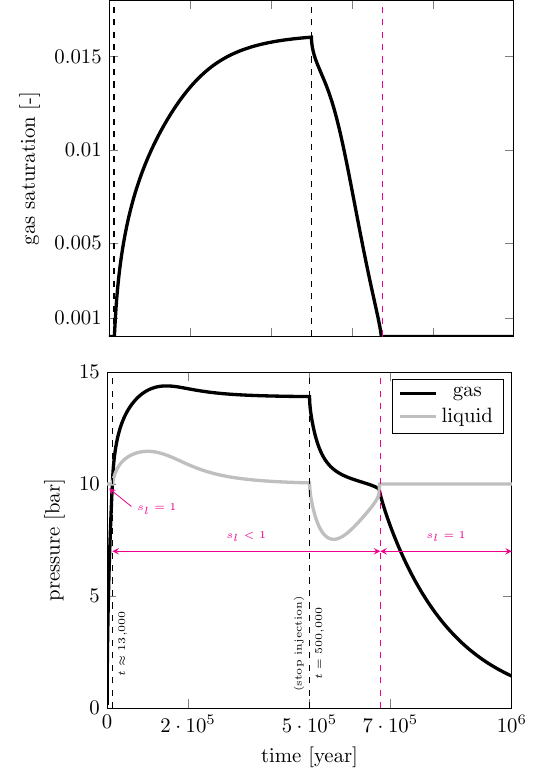}
\caption{Gas saturation (top) and liquid and gas pressures (bottom)
at the inflow boundary.}
\label{fig:inflow}
\end{figure}

\begin{figure}[ht]\centering
  \includegraphics[width=0.4\textwidth]{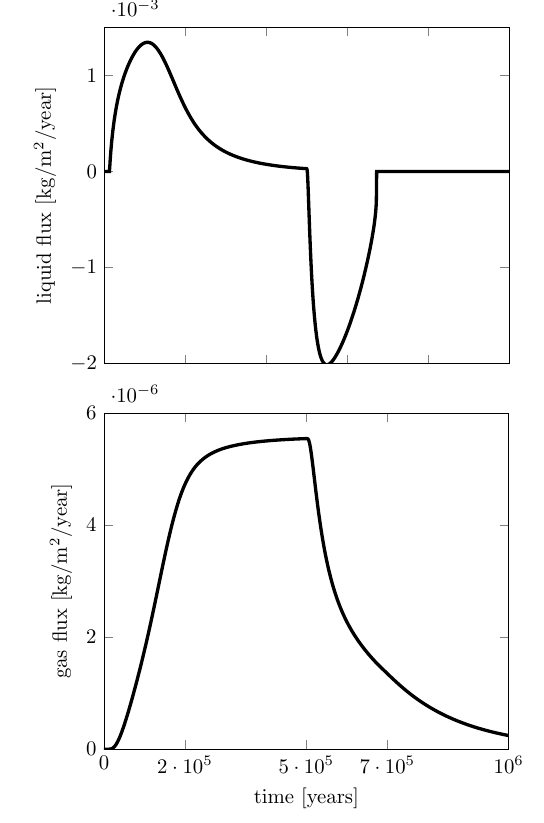}
  \caption{Liquid (top) and gas (bottom) flux at the outflow boundary.}
  \label{fig:bc_data}
\end{figure}

To validate our solver, we run simulations with the nominal parameters
and report the phase pressures and gas saturation at the inflow
boundary. Our results are consistent with those
reported in~\cite{Angelini2011,Bourgeat2009,Neumann2012}.
Figure~\ref{fig:inflow} shows the gas saturation (top) and the phase
pressures (bottom), with respect to time (years) during and after
injection. For $0 < t <13\times10^3$ years, the gas saturation
remains zero, all injected hydrogen dissolves into the liquid phase, the
whole domain is saturated with water, and the liquid pressure remains
constant.  At $t \approx 13\times10^3$ years, the maximum
solubility is reached and the gas phase appears at the injection
boundary. Gas saturation keeps growing during the period of hydrogen injection.
When injection stops at $t = 5\times 10^5$ years, gas
saturation decreases until it disappears. A negative water
flux is observed (see Figure~\ref{fig:bc_data}) as water comes from right to
left to fill in the empty space. At the end of the simulation, the
gas pressure continues to decrease and the liquid pressure gradient goes to
zero, as the system reaches a steady state.
%
% --------------------------  figure  ------------------------------------
%
% --------------------------  figure  ------------------------------------

\section{Modeling under uncertainty}\label{sec:model_uncertainty}
We seek to understand the impact of uncertainty in heterogeneous material
properties on model predictions. Specifically, we focus on uncertainties in
porosity and absolute permeability.  Our goal is to understand the impact of
uncertainties in material properties on the gas phase appearance/disappearance
in a two phase flow produced by hydrogen injection through a porous medium,
which is initially fully saturated with water.

\subsection{Modeling uncertainty in material properties}\label{sec:uncertainty_parameterization}
While in the setup of the benchmark problem constant values for porosity and
permeability were used, allowing for spatially varying porosity and
permeability provides a more realistic representation. This leads
to representation of these quantities as random fields.

\begin{figure}[ht]
  \centering
  \includegraphics[width=.45\textwidth]{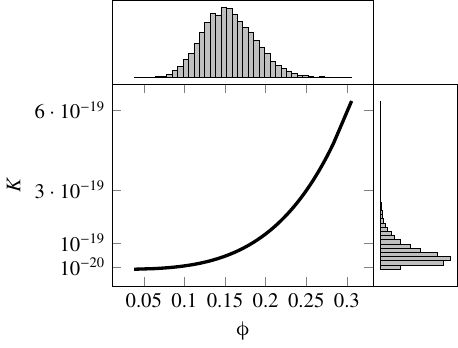}
  \includegraphics[width=.4\textwidth]{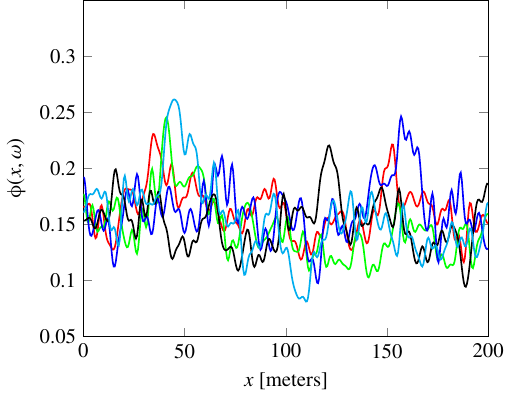}
  \caption{Top: the porosity permeability relation and the distributions of
  pointwise porosity and permeability. Bottom: a few realizations of the porosity field.}
  \label{fig:por_dist}
\end{figure}

We model the porosity, $\phii$, as a random field as follows.  Let $Z(x, \omega)$
be a Gaussian process, with exponential covariance function
$c(x, y) = e^{-|x - y|/\ell}$, where $\ell > 0$ is the correlation length. We
chose $\ell = 10$ \si{\meter} (recall the length of the domain is $200$
\si{\meter}). The covariance operator of $Z$ is defined by
\begin{equation}\label{equ:Cpar}
[\Cpar u](x) = \int_\Omega c(x, y) u(y) \,dy, \quad u \in L^2(\Omega).
\end{equation}
We define the random porosity field by
\begin{equation}\label{equ:por}
   \phii(x, \omega) = F_B^{-1} \left(F_G(Z(x, \omega)); \alpha_\text{beta}, \beta_\text{beta}\right).
\end{equation}
Here $F_B^{-1}(\cdot; \alpha_\text{beta}, \beta_\text{beta})$ is the inverse CDF of a 
Beta$(\alpha_\text{beta},
\beta_\text{beta})$ distribution and $F_G$ is the CDF of a standard normal distribution.
This ensures that for every $x \in \Omega$ the porosity is distributed according to
Beta$(\alpha_\text{beta}, \beta_\text{beta})$. The random permeability field is obtained using
a Kozeny--Carman relation~\cite{Costa06,Lie19}:
\begin{equation*}
  K(\phii) \propto \frac{\phii^3}{(1-\phii)^2}.
\end{equation*}
We set the proportionality constant in the above relation so that
$K(\bar\phii)= \bar K$, where $\bar\phii$ and $\bar K$ are the
nominal porosity and permeability values listed in Table~\ref{tab:nomvals}~(left). The
values of $\alpha_\text{beta}$ and $\beta_\text{beta}$ in~\eqref{equ:por} are set
such that the mode of the porosity distribution (at each $x \in \Omega$)
is the nominal porosity of $\bar\phii = 0.15$. Specifically, we chose
$\alpha_\text{beta} = 20$ and found $\beta_\text{beta}$ from the formula for the mode
of a Beta distribution:
$(\alpha_\text{beta} - 1)/(\alpha_\text{beta} +\beta_\text{beta} - 2) = \bar\phii$.
We depict the distributions for pointwise porosity and permeability values
along with the porosity permeability relation in Figure~\ref{fig:por_dist}~(top).
We note that the present setup provides a physically meaningful range of values
for porosity and permeability, for the application problem under study.
%
% --------------------------  figure  ------------------------------------

To facilitate uncertainty quantification, we consider a truncated KLE
of the Gaussian random field $Z(x, \omega)$ used in definition of $\phii(x,
\omega)$ in~\eqref{equ:por}. That is, we consider
\begin{equation}\label{equ:input_KLE}
  Z(x, \omega) \approx \sum_{i = 1}^\Np \sqrt{\lambda}_i \xi_i e_i(x),
\end{equation}
where $(\lambda_i, e_i)$, $i = 1, \ldots, \Np$ are the eigenpairs of the
covariance operator of $Z(x, \omega)$; see
e.g.,~\cite{AlexanderianReeseSmithEtAl18,BetzPapaioannouStraub14,Maitre10} for
details about the use of KL expansions for representing random fields in
mathematical models. For the present problem, we let $\Np = 100$,
which enables capturing over $96$ percent of the average variance of the process.
Notice that with the present setup, the uncertainty in the porosity
field is fully captured by the vector $\vec \xi =
[\begin{matrix} \xi_1 & \xi_2 & \cdots & \xi_\Np\end{matrix}]^\top$,
where $\xi_i$'s are the KLE coefficients in~\eqref{equ:input_KLE}, 
which are independent standard normal random variables.
As an illustration, we show a few realizations of the random porosity field in
Figure~\ref{fig:por_dist}~(bottom).

\subsection{The QoIs under study}
\begin{figure}\centering
\begin{tabular}{r}
\includegraphics[width=0.35\textwidth]{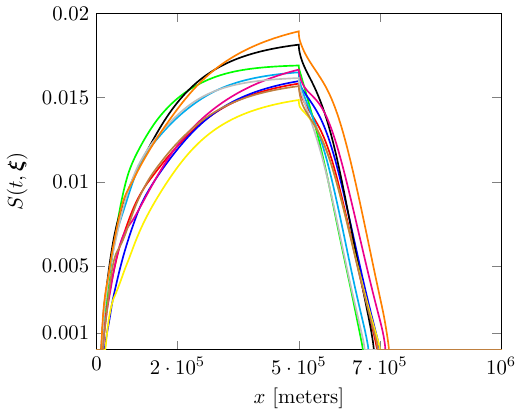}
\\
\includegraphics[width=0.37\textwidth]{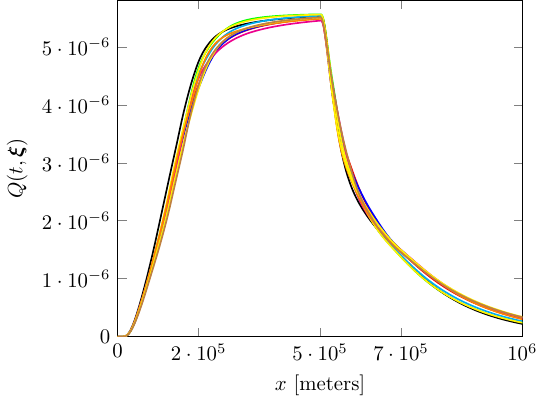}
\end{tabular}
\caption{A few realizations of the time evolution of top: gas saturation at the
inflow boundary, bottom: gas flux at the outflow boundary.}
\label{fig:QoIs}
\end{figure}

We focus on dynamics of hydrogen in gas phase by considering on the time
evolution of gas saturation and pressure at the inflow boundary and gas
flux at the outflow boundary. The units for gas pressure and gas flux are
[\si{bar}] and [\si{\kg\per\meter^2\per{year}}], respectively. These
time-dependent QoIs are indeed random field quantities due to randomness in
porosity and permeability fields.  Notice that since the uncertainty in
porosity field is encoded in the coefficients $\vec{\xi}$
in~\eqref{equ:input_KLE}, the randomness in these QoIs is also parameterized by
the vector $\vec\xi$ of the KL coefficients.  We denote the uncertain gas
saturation at the inflow boundary and gas flux at the outflow boundary by
$S(t, \vec\xi)$, and $Q(t, \vec\xi)$, respectively. In Figure~\ref{fig:QoIs},
we depict a few realizations of these uncertain QoIs.

\begin{figure}\centering
\includegraphics[width=0.35\textwidth]{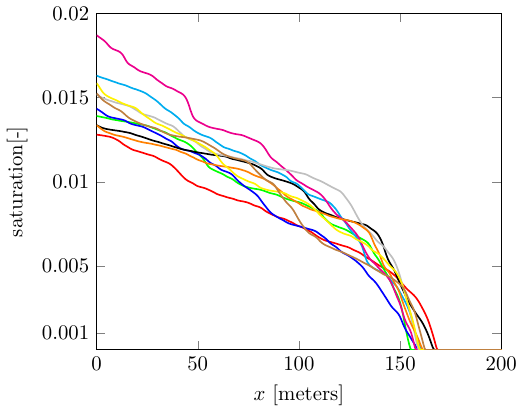}
\caption{Gas saturation at $t^* =300{,}091$ years.}
\label{fig:gas_sat_Tij}
\end{figure}

\begin{figure}\centering
\includegraphics[width=0.4\textwidth]{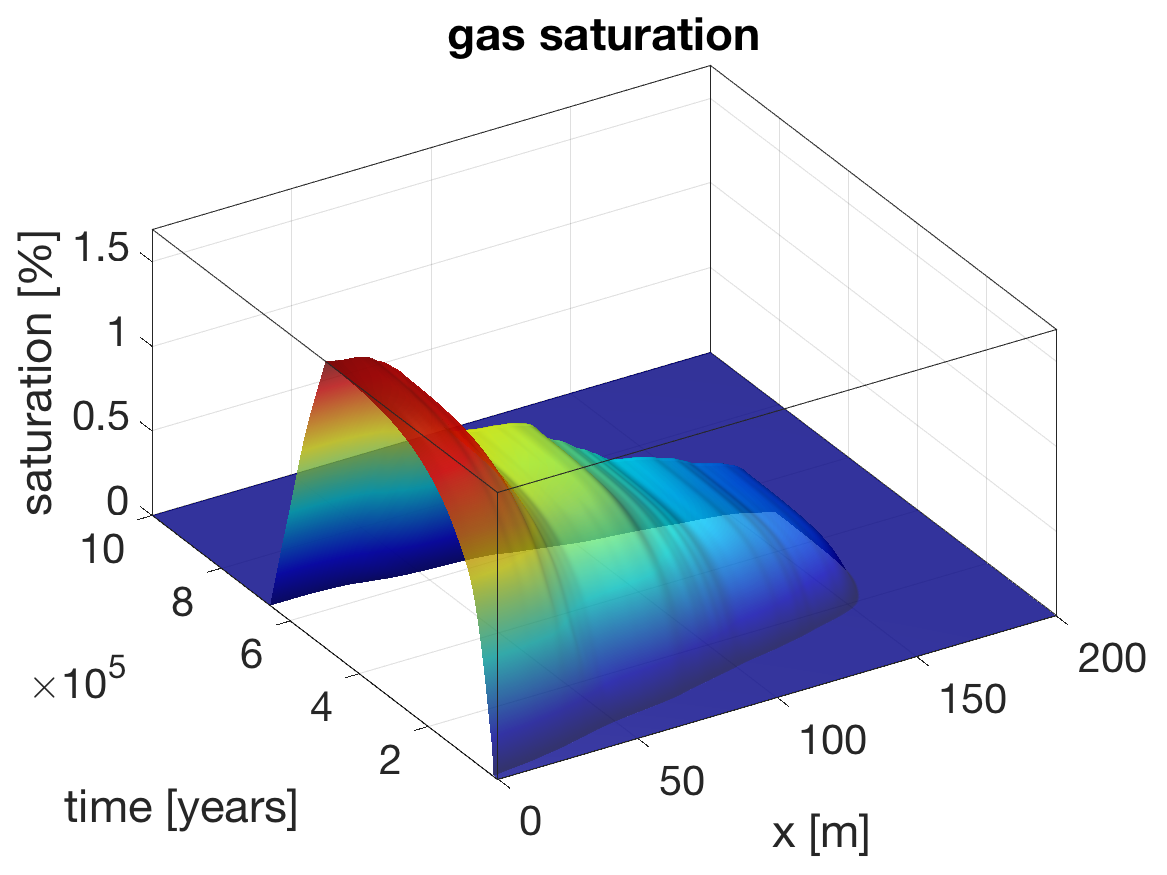}
\caption{Space time evolution of gas saturation.}
\label{fig:satg_surf}
\end{figure}

We also consider the gas saturation throughout the domain, at various points in
time. We denote this QoI by $S(x, \vec\xi; t^*)$, where
$t^*$ is a fixed time. Figure~\ref{fig:gas_sat_Tij} shows
a few realizations of this QoI at $t^* =  300{,}091$ years.  To further
illustrate the impact of spatial heterogeneity on the flow model, we also
report a plot of the gas saturation in  the space-time domain in
Figure~\ref{fig:satg_surf}.

Performing statistical studies and predictions on the QoIs outlined above is
challenging due to the high cost of solving the governing equations and the
high-dimensionality of the input and output spaces. A major aim of this article
is to present a surrogate modelling framework that approximates the
time- or space-dependent QoIs efficiently by reducing the input and output dimensions 
and using suitable approximations. 

% --------------------------  figure  ------------------------------------
% --------------------------  figure  ------------------------------------

\section{Spectral representations of random processes}
\label{sec:spectral}
\subsection{Karhunen Lo\'{e}ve expansions}\label{sec:KLE}
Here we discuss spectral representations of a function-valued output $f(s,
\vec\xi)$. We assume $f$ is a mean-square continuous random process. Such
processes admit spectral representations, as given by a Karhunen Lo\'{e}ve
expansion~(KLE)~\cite{Loeve77,Maitre10}:
\begin{equation}\label{eq:KLE_generic}
  f(s, \vec{\xi}) = \bar{f}(s)
  + \sum_{i=1}^{\infty} \sqrt{\lambda_i} f_i(\vec{\xi}) \Phi_i(s).
\end{equation}
Here $\bar{f}(s)$ is the mean of the process, $(\lambda_i, \Phi_i)$ are
the eigenpairs of the covariance operator $\Cqoi$ of the process,
\begin{equation}\label{equ:eigenvalue_problem}
    \Cqoi \Phi_i = \lambda_i \Phi_i, \quad i = 1, 2, \ldots,
\end{equation}
and
$f_i(\vec{\xi})$ are the KL modes,
\begin{equation*}
  f_i(\vec{\xi}) = \frac{1}{\sqrt{\lambda_i}} \int_\X \big(f(s, \vec{\xi})
  -\bar{f}(s)\big)\Phi_i(s)\, ds, \quad i = 1, 2, 3, \ldots
\end{equation*}
An approximation $\fn(s, \vec\xi)$ to $f(s, \vec\xi)$ can be obtained by
truncating~\eqref{eq:KLE_generic} and retaining the first $\Nqoi$ terms in the
series. In many physical and biological models the eigenvalues of $\Cqoi$
decay rapidly. Consequently, such QoIs can be represented with sufficient
accuracy by a truncated KLE with a small $\Nqoi$. Such processes are referred
to as ``low-rank''.

We rely on  Nystr\"{o}m's method to compute the KLE~\cite{Kress14}.
This approach, as used in the present work, requires sample averaging to
approximate the covariance kernel, because we do not in general have a
closed-form expression for the output covariance operator.  Typically, a
modest number of QoI evaluations is sufficient for accurately estimating the
dominant eigenpairs of the covariance operator $\Cqoi$. To determine a suitable
value for the number $\Nqoi$ of terms in a truncated KLE, we consider
\begin{equation}\label{eq:rk}
  r_k = \frac{\sum_{i=1}^k \lambda_i}{\sum_{i=1}^\infty \lambda_i}.
\end{equation}
The quantity $r_k$ represents the fraction of the average variance of $f$
captured by the first $k$ eigenvalues.
The steps for computing the truncated KLE of $f$ are included
in Algorithm~\ref{alg:KLE}, which is adapted
from~\cite{AlexanderianReeseSmithEtAl18}.

Note that evaluating the truncated KLE of $f$ requires computing the KL modes,
which in turn requires a model evaluation. To convert the truncated KLE into an
efficient surrogate model for $f$, we need a cheap-to-evaluate representation
for the KL modes. This approach is similar to the one taken
by~\cite{AlexanderianGremaudSmith20,LiIskandaraniLeHenaff16}, in which PCE
surrogates  are constructed for the modes of the related spectral
representations.  In section~\ref{sec:methods}, we modify this approach by
first reducing the dimension of the input parameter and then constructing the
KL modes surrogates in the reduced uncertain parameter space.

\subsection{Polynomial chaos expansions for $f_i(\vec\xi)$.}\label{sec:PCE_basic}
Recall, the polynomial chaos expansion of a square integrable function
$g(\vec{\xi})$ is a series approximation of the form
\begin{equation}\label{eq:generic_PC}
g(\vec{\xi}) \approx \sum_{k=0}^{\Npc} c_k{\Psi}_k(\vec{\xi}),
\end{equation}
where $\{{\Psi}_k\}_{k=0}^{\Npc}$ are a predetermined set of orthogonal
polynomials, and $\{c_k\}_{k=0}^{\Npc}$ are the corresponding expansion
coefficients~\cite{Maitre10}. Following a total order truncation~\cite{Maitre10}, $\Npc$ is given by
\begin{equation*}
\Npc + 1 = \frac{(\Nord + \Np)!}{\Nord!\Np!},
\end{equation*}
where $\Nord$ is the maximum total polynomial degree and $\Np$ is the dimension
of $\vec{\xi}$. There are a variety of approaches for determining the expansion
coefficients $\{c_k\}_{k=0}^{\Npc}$ including quadrature or regression based
methods~\cite{Maitre10}. For this application, we implement
sparse linear regression~\cite{Fajraoui17,Yan12}. In this method, the
expansion coefficients are found by solving
\begin{equation}\label{eq:sparse_lin_reg}
  \min_{\vec{c}\in\R^{\Npc}} \norm{\vec{\Lambda} \vec{c} -\vec{d}}_2^2, \quad
\text{subject to } \sum_{k=0}^{\Npc}|c_k|\leq \tau,
\end{equation}
where $\vec{\Lambda}\in \R^{\Ns \times \Npc}$ is defined by $\Lambda_{ij}
= {\Psi}_j(\vec{\xi}_i)$,\\ \noindent $\vec{d}
= (g(\vec{\xi}_1), g(\vec{\xi}_2), \dots, g(\vec{\xi}_{\Ns}))^\top$ is a vector
containing model evaluations, and $\tau$ is the sparsity control parameter.
Determining $\Nord$ and $\tau$ may be done with trial and error or with a
cross-validation process, as detailed in
Section~\ref{sec:numerics}.

\subsection{Bispectral surrogates}
Earlier we broached the subject of utilizing
PCEs to convert a truncated KLE into a surrogate model for $f$.
Consider the truncated KLE of $f$,
\begin{equation}\label{eq:trunc_KLE}
  \fn(s, \vec\xi) = \bar{f}(s) + \sum_{i=1}^{\Nqoi} \sqrt{\lambda_i}f_i(\vec\xi)
  \Phi_i(s).
\end{equation}
By replacing the KL modes in~\eqref{eq:trunc_KLE} with PCEs we construct a
surrogate model for $f$ of the form
\begin{equation}\label{eq:HD_surrogate}
\fPC(s, \vec\xi) = \bar{f}(s) + \sum_{i=1}^{\Nqoi} \sqrt{\lambda_i}
  f_i^{\PC}(\vec\xi)\Phi_i(s),
\end{equation}
where $f_i^{\PC}(\vec\xi)$ is the PCE for $f_i(\vec\xi)$, $i=1,\ldots, \Nqoi$.
Once constructed, a bispectral surrogate can be used to characterize
the statistical properties of the field QoI very efficiently.

To provide further insight, we also consider the  approximation error for a
bispectral surrogate.  Let $\norm{\cdot}$, represent the $L^2$ norm in the
product space $\paramspace\times\X$. The total error can be bounded as follows: 
\begin{align*}
   \norm{f-\fPC}^2 &\leq 2\norm{f-\fn}^2 + 2\norm{\fn - \fPC}^2\\
  &=2\sum_{i=\Nqoi+1}^\infty\hspace{-2mm}\lambda_i \\
  &+ \quad 2\sum_{i=1}^{\Nqoi} \lambda_i  \Big[\sum_{k=0}^{\Npc}(c_{i,k}-\hat{c}_{i,k})^2
\normLomega{\Psi_k}^2\Big]\\
& +\quad 2\sum_{i=1}^{\Nqoi} \lambda_i
  \Big[\sum_{j=1+\Npc}^{\infty}
  c_{i,j}^2\normLomega{\Psi_j}^2\Big].
\end{align*}
See Appendix~\ref{sec:append}, for a derivation of this bound. The first term
in the upper bound corresponds to KLE truncation error, the second term
corresponds to error due to inexact PCE coefficients, and the third term
corresponds to PCE truncation error.

Controlling the total error involves a balance between computational cost, accuracy
requirements, and the properties of the process. The KLE truncation error
gets smaller as $\Nqoi$ increases. However, increasing the number of terms in the
KLE increases the number of eigenpairs that need accurate approximations. Also, a
larger $\Nqoi$ results in more KL modes, each of which requires a
sufficiently accurate PCE.
Similarly, the PCE error can be minimized by increasing the maximum polynomial
degree, $\Nord$. However, this increases the total number of coefficients, which
increases the number of unknowns in~\eqref{eq:generic_PC}, resulting in
increased computational cost.

The function-valued QoIs in the present work are low-rank processes with a
high-dimensional input parameter. Therefore, a modest $\Nqoi$ will give a
sufficiently small KLE truncation error. However, for large $\Np$, estimating
the PCE coefficients for each KL mode with sufficient accuracy can become
computationally expensive. Our approach for addressing this challenge 
is presented in the next section.

%In the next section, we propose a method for alleviating the cost associated
%with PCE computation. First, we develop an efficient, parameter screening
%metric. Then, instead of constructing the PCEs for the KL modes in the full
%parameter space, we use our screening method to decrease the input parameter
%dimension and construct the KL mode PCEs on the resulting reduced parameter
%space.

\begin{algorithm}
  \caption{Computing the truncated KLE of $f$}.
  \label{alg:KLE}
  \begin{algorithmic}[1]
    \REQUIRE Quadrature nodes $s_k$ and
    weights $w_k$, $k = 1, \ldots, m$; Function evaluations
    $y_k^j = f(s_k, \vec{\xi}_j)$, $k = 1, \ldots, m$,
    $j = 1, \ldots, \Ns$; $r_k$ tolerance $0<tol<1$.
    \ENSURE Eigenpairs $\left(\lambda_i, \vec\Phi_i\right)$
    of the output covariance operator, 
    and KL modes evaluations $f_i(\vec{\xi_j})$, $j =1, \ldots, \Ns$,
    $i= 1, \ldots \Nqoi$.

    \STATE Compute mean $M_k =\frac{1}{\Ns} \sum_{j=1}^{\Ns}y_k^j$,
    $k = 1,\ldots m$.

    \STATE Center process $f_k^c(s_k,\vec{\xi}^j)=y_k^j - M_k$,
    $k = 1,\ldots, m$.

    \STATE Compute covariance matrix $\mat{C}$.\\
    $\mat{C}_{kl} = \frac{1}{\Ns-1}
    \sum_{j=1}^{\Ns} f_k^c(s_k, \vec{\xi}_j)f_l^c(s_l, \vec{\xi}_j)$,
    $k, l = 1, \ldots m$.

    \STATE Let $\mat{W} = \text{diag}(w_1,w_2, \ldots w_m)$ solve:\\
    $\mat{W}^{1/2}\mat{C}\mat{W}^{1/2} \vec{v}_k = \lambda_k\vec{v}_k$,
    $k=1, \ldots, m$.
    \STATE Determine $\Nqoi$.
    \FOR{$k=1,\ldots m$}
    \STATE Compute $r_k = \frac{\sum_{l=1}^k\lambda_l}{\sum_{l=1}^m \lambda_l}$.
    \IF{$r_k >tol$} \STATE $\Nqoi = k$; BREAK
    \ENDIF
    \ENDFOR
    \STATE Compute $\vec\Phi_k = \mat{W}^{-1/2}\vec{v}_k$, $k = 1, \ldots,
           \Nqoi$.
    \STATE Compute KL modes.\\
    $f_i(\vec{\xi}_j) = \frac{1}{\sqrt{\lambda_i}}
    \sum_{k=1}^{m} w_k f_k^c(s_k, \vec{\xi}_j) \vec\Phi_i(s_k)$,
    $i = 1, \ldots, \Nqoi$, $j =1, \ldots, \Ns$.
    \STATE Compute $\fn(s, \vec{\xi}_j) = \sum_{k=1}^{\Nqoi}
    \sqrt{\lambda_k} f_i(\vec{\xi}_j) \vec\Phi_k(s)$.
\end{algorithmic}
\end{algorithm}

\section{Method}
\label{sec:methods}
In this section, we present our approach for reducing the dimensionality of the
random vector $\vec{\xi} = [\begin{matrix}
\xi_1 & \xi_2 & \ldots & \xi_\Np\end{matrix}]^\top$  and
constructing a cheap-to-compute bispectral surrogate for
function-valued QoIs under study. We begin by describing a screening procedure for
input dimension reduction in Section~\ref{sec:screen}. Then, we discuss our
surrogate modeling approach that utilizes a truncated KLE of the output
(Section~\ref{sec:KLE}) along with generalized PCEs for the output KL
modes (Section~\ref{sec:PCE}). We also show how the surrogate model can be
used to efficiently compute the correlation function of the output, as well
as cross-correlation of two function-valued QoIs.

%
%
% SECTION
%
%
\subsection{Parameter screening}\label{sec:screen}
Consider a function-valued QoI $f(s, \vec{\xi}): \X \times \paramspace \to \R$,
where $\paramspace\subseteq \R^\Np$ is the sample space of the uncertain
parameters and $\X \subseteq \R^d$ is a compact set.  The set $\X$ can be
either a time interval, in which case $d = 1$, or a spatial region, in which
case $d \in \{1, 2, 3\}$. Here we consider the case of $d = 1$, as it applies
to our application problem, but the procedure below can be generalized to the
case of $d \in \{2, 3\}$ in a straightforward manner.

Parameter screening can be done using functional derivative-based global
sensitivity measures (DGSMs) given by \cite{CleavesAlexanderianGuyEtAl19}:
\begin{equation}\label{eq:func_DGSMs}
  \N_j(f) = \int_\X \int_\paramspace \left(\frac{\partial f(s, \vec\xi)}{\partial
  \xi_j}\right)^2 \, \mymu(d\vec\xi) \, ds, \quad j = 1, \ldots, \Np,
\end{equation}
where $\mymu$ is the law of the parameter vector $\vec\xi$.  These DGSMs can be
used to screen for ``unimportant'' inputs, which can be fixed at their
respective nominal values. These functional DGSMs, however, 
require gradient evaluations. For
complex models with high-dimensional parameters, such as the one 
considered
in the present work, gradient computation is challenging. While adjoint-based
gradient computation can overcome this, adjoint solvers are not always
available for complex flow solvers and implementing them may be infeasible.
Here we derive a screening indices based on ideas from active subspace
methods~\cite{Constantine15} and activity scores~\cite{ConstantineDiaz17} that
approximate the functional DGSMs and circumvent gradient computation.

Let us briefly recall the concept of the active subspace and activity
scores~\cite{ConstantineDiaz17}. Fix $s\in\X$ and let $(\lambda_k, \vec{u}_k)$,
$k = 1, \ldots, \Np$ be the eigenpairs of the matrix
\begin{equation}\label{eq:C}
  \mat{G}=\int_\paramspace [\nabla f(s, \vec\xi)][\nabla f(s, \vec\xi)]^\top \,
  \mymu(d\vec\xi),
\end{equation}
where we assume the eigenvalues are sorted in descending order.
In many cases there exists an $M$ such that $\lambda_M  \ll \lambda_{M+1}$,
representing a gap in the eigenvalues. The {\emph{active subspace}} corresponds
to the subspace spanned by eigenvectors $\{\vec{u}_k\}_{k=1}^M$; 
this subspace captures the directions in the uncertain parameter space along which
the QoI varies most. The case of a one-dimensional active subspace is
surprisingly common~\cite{Constantine15}.
The \emph{activity scores}~\cite{ConstantineDiaz17} utilize the active subspace
structure to provide approximate screening indices, given by
\begin{equation}\label{eq:actscores}
   \alpha_j[f(s, \cdot); M] =  \sum_{k=1}^M
   \lambda_k \ip{\vec{e}_j}{\vec{u}_k}^2,
   \quad j = 1, \ldots, \Np, \quad M \leq \Np,
\end{equation}
where $\ip{\cdot}{\cdot}$ denotes the Euclidean inner product
and $\vec{e}_j$ is the $j$th coordinate vector in $\R^\Np$. One can use the
activity scores to approximate functional DGSMs according to
\[
\int_\X \alpha_j[f(s, \cdot); M] \, ds, \quad j = 1, \ldots, \Np.
\]
Note that with $M = \Np$, we recover the exact DGSMs~\cite{ConstantineDiaz17}.
Computing the activity scores still requires gradient computation, as seen in
the definition of the matrix $\mat{G}$ in~\eqref{eq:C}.  For cases where full
model gradients are unavailable, we propose use of suitable and
cheap-to-compute surrogate models for the purpose of computing the activity
scores.  One possibility is the use of global linear models as done
in~\cite{Constantine15}, for the case of scalar QoIs.  Building on this idea,
we use a global linear model for $f(s, \vec\xi)$, use the gradient of the
linear model to approximate the matrix $\mat{G}$, and 
define screening indices for
function-valued QoIs. Specifically, we construct a global linear approximation
$\tilde{f}$ for the QoI
\begin{equation}\label{eq:lin_mod}
\tilde{f}(s, \vec\xi) = b_0(s) + \sum_{j = 1}^\Np b_j(s) \xi_j.
\end{equation}
Next, we use the activity scores for $\tilde{f}$ as a ``surrogate'' for the
scores of $f$. Note that $\nabla \tilde{f}(s, \vec\xi) = \vec{b}(s)$, where
$\vec{b}(s) = 
[\begin{matrix} b_1(t) & b_2(t) & \cdots & b_\Np(t)\end{matrix}]^\top$.
The matrix $\mat{G}(s)$ in~\eqref{eq:C}, using $\tilde{f}$ in place of $f$ then
simplifies to
$\mat{G}(s) = \vec{b}(s)\vec{b}(s)^\top$.
This rank one matrix can be written as
\begin{equation*}
  \mat{G}(s) = \lambda \vec{u}(s) \vec{u}(s)^\top,
\end{equation*}
where $\lambda = \norm{\vec{b}(s)}_2^2$, and
$\vec{u}(s) = \vec{b}(s)/\norm{\vec{b}(s)}_2$. (Here $\norm{\cdot}_2$
denotes the Euclidean vector norm.)  Hence, the corresponding active
subspace for $\tilde{f}$ is 1-dimensional resulting in activity scores
\begin{equation*}
  \tilde{\alpha}_j(s) = b_j^2(s),
  \quad j = 1, \ldots, \Np.
\end{equation*}
This gives rise to the following approximate functional DGSMs:
\[
  \tilde{\N}_j(f) := \int_\X  b^2_j(s) \, ds.
\]
This relationship motivates the following normalized screening indices
\begin{equation}\label{eq:screening_indc}
  \sj = \frac{\tilde{\N}_j(f)}{\sum_{l=1}^\Np \tilde{\N}_l(f)},
  \quad j = 1, \ldots, \Np.
\end{equation}
Henceforth, we refer to $\sj$ as the {\emph{screening index}} of $f$ with
respect to parameters $\xi_j$.

The purpose of the screening indices $\sj$ is to inform input parameter
dimension reduction. Let $K_r$ be an ordered index set with cardinality
$\np<\Np$, corresponding to parameters with a screening index above some
tolerance $tol \in (0, 1)$. We denote the reduced input parameter vectors
$\vec{\xi}^r$, where each component $\xi_i^r$, $i = 1, \ldots, \np$ corresponds
to the $i$th element of $K_r$.

Next, we discuss the computation of the global linear model for $f$.
This is done by computing a linear model
at each point $s\in \X$, which can be done efficiently using linear regression.
Recall that $\X$ is assumed to be a (compact) 
subset of $\R$ (i.e., in one space dimension). Specifically, 
we assume 
$\X = [s_0, s_F]$.
We discretize $\X$ using a grid
\begin{equation*}
   s_0 = s_1 < s_2 < s_3 < \cdots < s_m = s_F.
\end{equation*}
Denote $\bar{\vec{b}}(s) = 
[\begin{matrix} b_0(s) & b_1(s) & b_2(s) & \cdots
& b_{\Np}(s)\end{matrix}]^\top$, with $b_j$, $j = 0, \ldots, \Np$ as
in~\eqref{eq:lin_mod}. We require a set of model evaluations,
\begin{equation*}
  y_k^i = f(s_k, \vec\xi_{i}), \quad i = 1, \ldots, \Ns.
\end{equation*}
The number of samples required depends on computational budget as well as
the application problem under study. We show in our numerical results that a
modest $\Ns$ is adequate for the proposed approach, and the application
problem considered herein.

Let $\vec{y}_k = [\begin{matrix} y_k^1 & y_k^2 & \cdots
& y_k^{\Ns}\end{matrix}]^\top \in \R^{\Ns}$, and define the matrix
\begin{equation}
\label{eq:Amatrix}
  \mat{A} =
  \begin{bmatrix}
  1 & \vec\xi_1^\top \\
  1 & \vec\xi_2^\top \\
  1 & \vec\xi_3^\top \\
  \vdots & \vdots \\
    1 & \vec\xi_{\Ns}^\top
  \end{bmatrix}.
\end{equation}
The vectors $\bar{\vec{b}}(s_k)$ can be computed numerically by solving
linear least squares problems
\begin{equation}\label{eq:lin_least}
  \bar{\vec{b}}(s_k) = \argmin_{\vec{b} \in \R^{\Ns+1}}
  \norm{\mat{A} \vec{b} - \vec{y}_k}_2^2,
\end{equation}
for $k = 1, \ldots, m$. Note that here we assume $\mat{A}$ has full column rank and
we are in the overdetermined case, i.e. $\Ns>\Np+1$.
Under these assumptions, the QR factorization $\mat{A} =\mat{Q}\mat{R} $ may
be used to solve the linear regression problem in~\eqref{eq:lin_least} by
\begin{equation*}
 \bar{\vec{b}}(s_k) = \mat{R}^{-1} \mat{Q}^\top \vec{y}_k.
\end{equation*}
Then, for each $k = 1,\ldots, m$, the cost of computing $\vec{b}_k$ is one
matrix-vector product with $\mat{Q}^\top$ and one triangular solve.
The procedure for computing the global linear model is summarized in
Algorithm~\ref{alg:screening_ind_qr}. In the case where the dimension of
$\vec\xi$ is larger than the number of available function evaluations, i.e.
$\Ns< \Np+1$, other methods for solving the linear regression in
Equation~\eqref{eq:lin_least}, e.g. using SVD, can be used.

%%%%%%%%%%%%%%%%%%%%%%%%%%%%%%%%%%%%%%%%%%%%%%%%%%%%%%%
Note that when using a global linear model for the approximating function
$\tilde{f}$, the screening indices~\eqref{eq:screening_indc} coincide
with the normalized functional DGSMs of $\tilde{f}$,
\[
\sj = \frac{\mathbb{N}_j(\tilde{f})}{\sum_{i=1}^\Np \mathbb{N}_i(\tilde f)},
\]
with $\mathbb{N}_j$ is as in~\eqref{eq:func_DGSMs}.  
Additionally, since we
have independent standard normal input parameters the
screening indices are equal to the function-valued Sobol' indices
\cite{GamboaJanonKlein14,AlexanderianGremaudSmith20} 
of $\tilde{f}$, as well as
the square root of the standard regression
coefficients~\cite{Helton93,ConstantineDiaz17} for $f$.
In general, the relations across these sensitivity measures will not hold for
alternative choices of $\tilde{f}$ or input parameter distributions.  

We emphasize that, while simplifications occur when a global linear
model is used, the
proposed screening approach, i.e., computing activity scores 
of the $\tilde{f}$, is intended to be flexible and adjustable to
alternative modeling approaches for $\tilde{f}$. Furthermore, the proposed
screening method is not constrained by the assumption of independent parameters
and can be used for the case of dependent inputs.  The only requirements of
the proposed screening method are that $\tilde{f}$ be cheap
to compute and adequately approximate the full model for the purposes of
parameter screening.

\begin{algorithm}
\caption{Computing the screening indices $\sj$, $j = 1, \ldots, \Np$: the
  overdetermined case.}
  \label{alg:screening_ind_qr}
  \begin{algorithmic}[1]
\REQUIRE 
Quadrature nodes $s_k$ and weights $w_k$, $k = 1, \ldots, m$.
Function evaluations
$y_k^i = f(s_k, \vec\xi_{i})$, $i = 1, \ldots, \Ns$, $k = 1, \ldots, m$;
\ENSURE Sensitivity measures $\sj$, $j = 1, \ldots, \Np$.
\STATE Form the matrix $\mat{A}$ in~\eqref{eq:Amatrix} and compute its
QR factorization, $\mat{A} = \mat{Q}\mat{R}$.
\FOR{$k=1$ to $m$}
   \STATE Compute $\vec{z}_k = \mat{Q}^\top \vec{y}_k$.
   \STATE Solve $\mat{R}\bar{\vec{b}}(s_k) = \vec{z}_k$.
\ENDFOR
\FOR{$j=1$ to $\Np$}
    \STATE Compute $\tilde{\N}_j = \sum_{k=1}^\Np w_k b_j(s_k)^2$.
\ENDFOR
\FOR{$j=1$ to $\Np$}
    \STATE Compute $\sj = \tilde{\N}_j / (\sum_k \tilde{\N}_k)$.
\ENDFOR
\end{algorithmic}
\end{algorithm}

%
%
% SECTION
%
%
%
\subsection{Polynomial Chaos surrogates for KL modes}\label{sec:PCE}
To form a surrogate model, we construct a PC surrogate
$f_i^{\PC}(\vec{\xi}^r)$,$i = 1, \ldots \Nqoi$ in the \emph{reduced}
parameter space. Explicitly, we have the following training data for the
KL mode surrogates: the input parameter samples
$W = \{\vec{\xi}_j^r\}_{j=1}^{\Ns}$  and, for each KL mode $f_i$
$i = 1, \ldots \Nqoi$, we have the evaluations
$F_i = \{f_i(\vec{\xi}_j)\}_{j=1}^{\Ns}$. For each KL mode $f_i$, we use
the corresponding training data
to solve the optimization problem~\eqref{eq:sparse_lin_reg} for the
coefficients $\vec{c}$; see Algorithm~\ref{alg:KLE_surrogate} for more details.
Observe that each input parameter sample $\vec{\xi}_j^r$ is the reduced version
of the original input parameter sample, whereas the data points in $F_i$
correspond to the KL mode $f_i$ evaluated on the full parameter $\vec{\xi}_j$.
Utilizing the data this way has two benefits. Firstly, we do not require more
model evaluations. Secondly, the KL modes corresponding to the exact
QoI capture the behavior of $f$ more accurately than the KL
modes corresponding to an $f$ re-evaluated in the reduced parameter space. After
the PCE for each KL mode is computed, we replace each $f_i(\vec{\xi})$ in the KL
expansion~\eqref{eq:KLE_generic} with the corresponding $f_i^{\PC}$ to form a
(reduced space)
bispectral surrogate  for $f$:
\begin{equation}\label{eq:surrogate}
  f(t, \vec{\xi}) \approx \fPC(t, \vec{\xi}^r) = \hat{f}(t) +
  \sum_{i=1}^{\Nqoi} \sqrt{\lambda_i} f_i^{\PC}(\vec{\xi}^r) \Phi_i(t).
\end{equation}
In section~\ref{sec:numerics} we demonstrate the proposed approach for
dimension reduction and surrogate modeling for temporally varying QoI
$S(t, \vec{\xi})$ and $Q(t, \vec{\xi})$, as well as spatially varying QoI
$S(x, \vec \xi)$.

Bispectral surrogates of the form~\eqref{eq:surrogate} can be sampled efficiently to
study the statistical properties of the QoI. As seen below, such surrogates can also be
used to efficiently compute the correlation structure of function-valued
outputs.

%%%%%%%%%%%%%%%%%%%%%%%%%%%%%%%%%%%%%%%%%%%%%%%%%%%%%%%

Note that the use of a PCE to construct the bispectral surrogate for
$f$ relies on the assumption of independent input parameters.  In the case of
dependent inputs, an alternative surrogate modeling approach should be utilized
for the output KL modes. The outline of the procedure would remain the same,
with the final surrogate model being constructed in the reduced parameter
space.

%Add sentence about independence as it applies to the screening approach as as it applies to the surrogate modeling approach.

%BUT if not PCE, then the surrogate model is no longer bispectral ….

%Look, if the inputs are not independent, more generally, on could replace with more general surrogate models.

%But, it will not be bispectral anymore….

%Note that the bispectral surrogates just “spectral”(?) can be explained in conclusion?

%%%%%%%%%%%%%%%%%%%%%%%%%%%%%%%%%%%%%%%%%%%%%%%%%%%%%%%

%
\begin{algorithm}
\caption{Computing the surrogate model $\fPC$}.
  \label{alg:KLE_surrogate}
  \begin{algorithmic}[1]
    \REQUIRE Reduced input parameters $\vec{\xi}_j^r\in\R^{\np}$,
    $j = 1,\ldots \Ns$;
    KL mode evaluations
    $f_i^k = f_i(\vec{\xi}_j)$, $i = 1, \ldots, \Nqoi$, $j = 1, \ldots, \Ns$;
    highest polynomial degree $\Nord$; sparsity parameter $\tau$;
    polynomial basis $\Psi_k$, $k = 1\ldots \Npc$.
    \ENSURE Surrogate model $\fPC(t, \vec{\xi}^r)$ and polynomial KL mode
    expansions $f_i^{\PC}(\vec{\xi}^r)$, $i = 1,\ldots, \Nqoi$.
    \FOR{$i=1$ to $\Nqoi$}
    \STATE Let $\vec{d}_i = [f_i(\vec{\xi}_1), \ldots, f_i(\vec{\xi}_{\Ns})]$\\
    and $\mat{\Lambda}_{kj} = \Psi_{k}(\vec{\xi}_j^r)$.
    \STATE Solve
    $$\min_{\vec{c}_i\in\R^{\Npc}} \norm{\vec{\Lambda} \vec{c}_i -\vec{d}_i}_2^2,
    \qquad \text{subject to } \sum_{k=0}^{\Npc}|c_k|\leq \tau$$
    \ENDFOR
    \STATE Form $f_i(\vec{\xi}^r) =
    \sum_{k=1}^{\Npc} c_{i,k}\Psi_k(\vec{\xi}^r)$, $i = 1, \ldots, \Nqoi$.
    \STATE Form $\fPC(s, \vec{\xi}^r) = \sum_{k=1}^{\Nqoi}
    \sqrt{\lambda_i(\mat{C})}f_i^{\PC}(\vec{\xi}^r)\Phi_i(s)$.

\end{algorithmic}
\end{algorithm}

%
%
% SECTION
%
%
%
\subsection{Correlation structure of the output}\label{sec:correlation}

Let $f:\X\times\paramspace\to\R$ be a random process with mean
$\bar{f}(s)$ and  assume $f$  admits a surrogate $\fPC$ of the form
in~\eqref{eq:surrogate}. 
It is straightforward to show that the covariance operator of 
$\fPC$ satisfies 
\begin{align}
  c_{f}(s_1, s_2) &=  \Cov{\fPC(s_1, \cdot), \fPC(s_2, \cdot)}\nonumber \\
  &=\sum_{i=1}^{\Nqoi} \sum_{j=1}^{\Nqoi}
  \sum_{k=1}^{\Npc} \coma_i^k\coma_j^k \norm{\Psi_k}_{L^2(\paramspace)}^2 
   \Phi_i(s_1)\Phi_j(s_2)
  \label{eq:cov_1},
\end{align}
for $\coma_i^k =c_{i,k}\sqrt{\lambda_i}$ and $\norm{\cdot}_{L^2(\paramspace)}$ 
denotes the $L^2$
norm on $\paramspace$. Let us define
\[
  B_{ij} := \sum_{k=1}^{m} \coma_i^k\coma_j^k 
            \norm{\Psi_k}_{L^2(\paramspace)}^2, \quad
i, j = 1, \ldots, \Nqoi,
\]
and
\begin{align*}
  \vec{p}(s) &:= [\begin{matrix} \Phi_1(s) ~ \Phi_2(s)~\ldots
~\Phi_{\Nqoi}(s) \end{matrix}]^\top.\\
\end{align*}
We can rewrite the expression in~\eqref{eq:cov_1} as
\begin{equation*}
  c_{f}(s_1, s_2) =
\ip{\vec{p}(s_1)}{\mat{B}\vec{p}(s_2)}, \label{eq:cov_analytic}
\end{equation*}
where $\ip{\cdot}{\cdot}$ denotes the Euclidean inner product.
Using this, we an also obtain the \emph{correlation function} of $\fPC$:
\begin{equation}\label{eq:corr}
  \rho_f(s_1, s_2) = \frac{c_f(s_1, s_2)}{\sqrt{c_f(s_1, s_1)}
  \sqrt{c_f(s_2, s_2)}}.
\end{equation}

We can also compute the cross-covariance
function of two random processes represented
via bispectral surrogates. Consider
a random process $g$ approximated by the surrogate model
\begin{equation*}
  \gPC = \bar{g}(s) + \sum_{j=1}^{\Mqoi}\sum_{k=0}^{\Mpc}\sqrt{\gamma_j}
  d_{j,k}\Psi_k(\vec\xi^r) \tilde{\Phi}_j(s),
\end{equation*}
where $\Mqoi$ is the number of KL modes, $(\gamma_j, \tilde{\Phi}_j(s))$ are
the eigenpairs corresponding to the covariance function of $g$, $\Mpc$ is the
maximum polynomial degree, and $d_{i,k}$ are the PCE coefficients.
A calculation similar to the one above gives the cross--covariance function of
$\fPC$ and $\gPC$ as
\begin{equation*}
  c_{fg}(s_1, s_2)
  =\ip{\vec{p}(s_1)}{\tilde{\mat{B}}\vec{q}(s_2)}, \label{eq:cross_analytic}
\end{equation*}
where
\[
\begin{aligned}
\vec{q}(s) &:= [\begin{matrix} \tilde{\Phi}_1(s)~\tilde{\Phi}_2(s)
~\ldots~\tilde{\Phi}_{\Mqoi}(s) \end{matrix}]^\top,\\
  \tilde{B}_{i,j} &:= \sum_{k=1}^{m} \coma_i^k\comb_j^k 
  \norm{\Psi_k}_{L^2(\paramspace)}^2,
  \quad i = 1, \ldots, \Nqoi,~j = 1, \ldots \Mqoi,
\end{aligned}
\]
with $\comb_j^k = d_{j,k}\sqrt{\gamma_j}$. We can also compute the
cross-correlation function,
\begin{equation}\label{eq:crosscorr}
  \rho_{fg}(s_1, s_2) =
  \frac{c_{fg}(s_1, s_2)}{\sqrt{c_f(s_1,s_1)} \sqrt{c_g(s_2, s_2)}},
\end{equation}
where $c_g$ is the covariance function of $\gPC$.

\section{Numerical Results}
\label{sec:numerics}
In this section, we demonstrate the dimension reduction and surrogate modeling
approach proposed in Section~\ref{sec:methods} for temporally and spatially
varying QoIs discussed in Section~\ref{sec:model_uncertainty}. In
Section~\ref{sec:gas_saturation_surr}, we detail surrogate model construction
for gas saturation at the inflow boundary.  To provide further insight, we also
consider surrogate modeling for gas flux at the outflow boundary in
Section~\ref{sec:gas_flux_surr} and for gas saturation across the spatial
domain in  Section~\ref{sec:spatial_surr}.  Finally, in Section~\ref{sec:use_surrogate},
we use the surrogates constructed in Section~\ref{sec:gas_saturation_surr} and
\ref{sec:gas_flux_surr} to better understand the behavior and properties of the
corresponding QoIs.
%
%                             SECTION
%
\subsection{Gas saturation at the outflow boundary}
\label{sec:gas_saturation_surr}
Here we focus on gas saturation at the inflow boundary, i.e., $S(t,
\vec{\xi})$.  Recall that the input parameter $\vec\xi$ parameterizes the
uncertainty in the porosity field, as described in
Section~\ref{sec:uncertainty_parameterization}, and has dimension $\Np=100$.
For the present numerical study, we computed a database of $550$ model
evaluations, which we use for parameter screening and surrogate model
construction.

\begin{figure}[h!]
\includegraphics[width=0.4\textwidth]{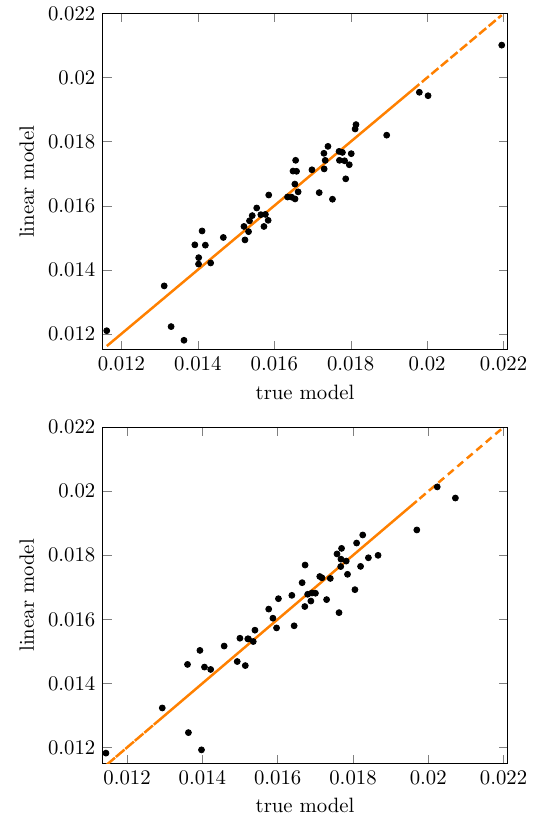}
\caption{Fifty point comparison of the true model to the linear model for
$S(t, \vec\xi)$  at top: $t=400{,}234$ years, bottom: $t=500{,}106$ years.}
\label{fig:linear_model_gas_saturation}
\end{figure}

\begin{figure}[h!]
\includegraphics[width=0.4\textwidth]{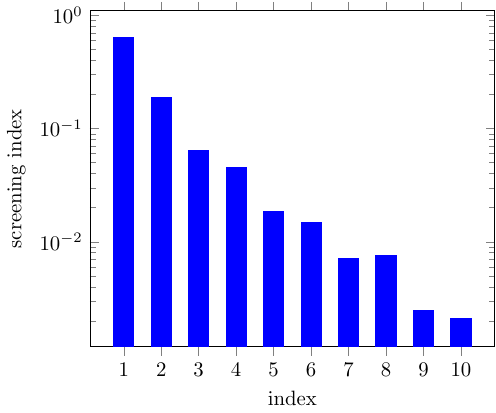}
\caption{Screening indices $\sj$ for $S(t, \vec\xi)$ calculated using
Algorithm~\ref{alg:screening_ind_qr} with $500$ full QoI samples. Indices above
$tol = 0.002$ displayed only.} \label{fig:screening_gas_saturation}
\end{figure}

%                             SECTION
%
\textbf{Input parameter screening.}
We use Algorithm~\ref{alg:screening_ind_qr} with composite trapezoid rule and $\Ns = 500$ full model
evaluations to compute the screening indices $\sj$, $j = 1, 2, \ldots, \Np$, for
$S(t, \vec\xi)$.  The remaining $50$ realizations were used for validation of
the linear models computed as a part of the algorithm.  In
Figure~\ref{fig:linear_model_gas_saturation}, we report representative
comparisons of the linear model versus the exact model, at the validation
points at selected times. Note that the linear models capture
the overall behavior of the model response.
%
%------------------------------- fig ----------------------------------------
In Figure~\ref{fig:screening_gas_saturation}, we report the screening indices
that are above the importance threshold $tol = 0.002$. The parameters with
screening indices below $tol$ are considered unimportant.  This reduces
the input parameter dimension from $\Np=100$ to $\np = 10$ and the resulting
reduced parameter is $\vec\xi^r = [\begin{matrix}
\xi_1 & \ldots & \xi_{10}\end{matrix}]^T$.
%------------------------------- fig ----------------------------------------

%                             SECTION

\begin{figure}[h!]
\includegraphics[width=0.45\textwidth]{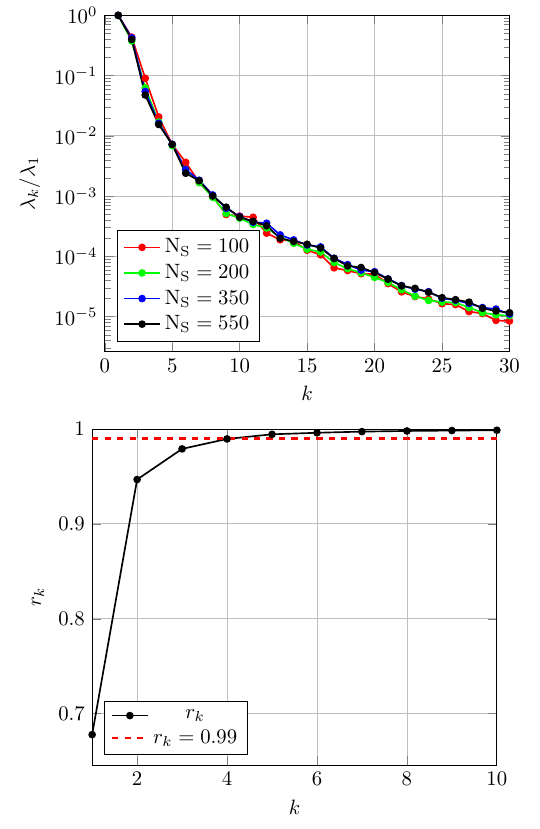}
%------------------------------- fig ----------------------------------------
\caption{Top: comparison of ratio $\lambda_k/\lambda_1$, $k = 1\ldots, 30$ for
$\lambda_i(\Cqoi)$ corresponding to $S(t, \vec\xi)$ computed with various sample
sizes, bottom: $r_k$ as defined in~\eqref{eq:rk}, $k = 1,\dots, 10$, for
$S(t, \vec\xi)$. Dotted line corresponds to $0.99$.}
\label{fig:eigs_gas_saturation}
\end{figure}

\textbf{Spectral representation of the QoI}.
Next, we compute the KLE of $S(t, \vec\xi)$ using
Algorithm~\ref{alg:KLE}. This requires solving the
eigenvalue problem~\eqref{equ:eigenvalue_problem}, with $\Cqoi$ being the
covariance operator of $S(t, \vec\xi)$.  We use a sample average approximation
of $\Cqoi$ with sample size $\Ns \in \{100, 200, 350, 550\}$ exact QoI
evaluations, as detailed in Algorithm~\ref{alg:KLE}. For this calculation we 
utilize the weights and nodes associated with the composite trapezoid rule. In
Figure~\ref{fig:eigs_gas_saturation}~(top), we show the computed (dominant)
eigenvalues of $\Cqoi$.  We note that the dominant eigenvalues are
approximated well even with $\Ns = 100$. We use the computations corresponding
$\Ns = 550$ in what follows. We note that the eigenvalues of the output covariance
operator decay rapidly. We also report $r_k$ from equation~\eqref{eq:rk},
in Figure~\ref{fig:eigs_gas_saturation}~(bottom). We note that $r_k$ exceeds $0.99$
for $k \geq 5$. This indicates that $S(t, \vec\xi)$ is a low-rank process and a KL
expansion with $\Nqoi = 5$ provides a suitable approximation of the QoI.
Consequently, we consider the truncated KL expansion of $S(t, \vec\xi)$
%------------------------------- fig ----------------------------------------
\begin{equation}\label{eq:gas_saturation_KLE}
\Sn(t, \vec{\xi}) = \bar{S}(t) + \sum_{i=1}^{\Nqoi}
  \sqrt{\lambda_i(\Cqoi)} S_i(\vec{\xi}) \Phi_i(t),
\end{equation}
where $\Nqoi = 5$. The next step is to compute
PCEs for the KL modes $S_i(\vec\xi)$, $i = 1, \ldots, \Nqoi$.
%
%                               SECTION

\textbf{PCE surrogates of the KL modes}.
Next, we construct a bispectral surrogate for $S(t, \vec\xi)$ which we
denote $\SPC$.
Recall that the components of $\vec{\xi}^r$ are sampled from a Gaussian
distribution.  Hence, we utilize the $\np$-variate Hermite polynomials as the
orthogonal basis for the PC expansions, with $\np = 10$.  We use the
sparse-regression approach (see Section~\ref{sec:PCE_basic}) for computing PCEs
of the output KL modes (see Section~\ref{sec:PCE}).  To determine suitable
values for the maximum polynomial degree $\Nord$ and the sparsity parameter
$\tau$, we use a 10-fold cross validation procedure, which we briefly explain next.

%
%------------------------------- fig ----------------------------------------
\begin{figure}[h!]
\includegraphics[width=0.45\textwidth]{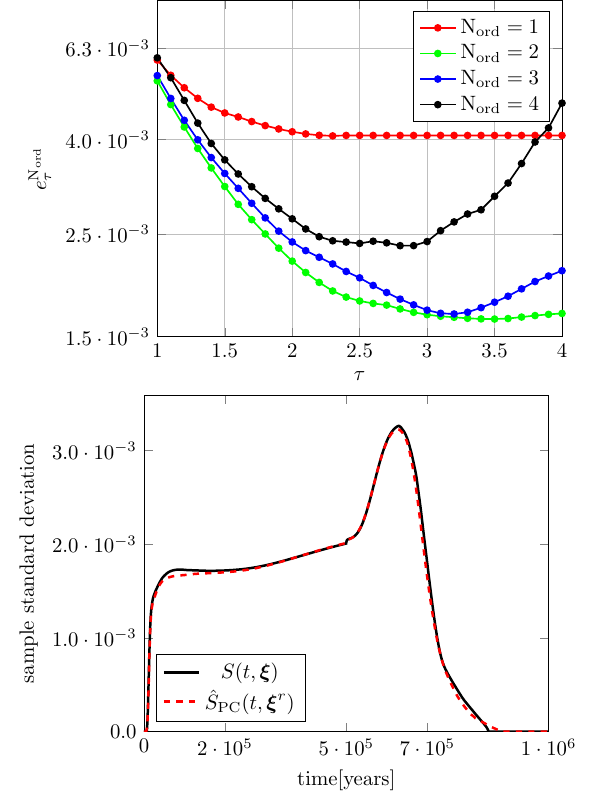}
%------------------------------- fig ----------------------------------------
\caption{Top: cross-validation results for $\tau = \{1, 1.1, \ldots, 3.9, 4\}$ and
$\Nord = \{1, 2, 3, 4\}$ for gas saturation, bottom: comparison of sample standard
deviations of $S(t, \vec\xi)$ and $\SPC(t, \vec\xi^r)$ computed on $200$ sample
points.}
\label{fig:cross_validation_gas_saturation}
\end{figure}
Note that for each evaluation of $\Sn(t, \vec\xi_j)$, $j = 1, \ldots, \Ns$, there
is a corresponding KL mode evaluation $S_i(\vec\xi_j)$, for $i=1, \ldots, \Nqoi$.
We separate the parameter samples into $W = \{\vec\xi_j^r\}_{j=1}^{350}$ and
$\hat{W} = \{\vec\xi_j^r\}_{j=351}^{550}$. Similarly, for each
$i = 1, \ldots \Nqoi$, we have $F_i = \{S_i(\vec\xi_j)\}_{j=1}^{350}$ and
$\hat{F}_i =\{S_i(\vec\xi_j)\}_{j=351}^{550}$.

We partition $W$ and $F_i$, $i =1, \ldots, \Nqoi$ $10$ different ways, such
that each data partition has a $35$ point validation set and a $315$ point
training set. Let $W^k$ and $F_i^k$, denote the $k$th such data partition, $k =
1, \ldots, 10$.  Next, for every combination of $\Nord \in \{1,\dots 4\}$,
$\tau \in \{1, 1.1, 1.2, \ldots, 3.9, 4\}$, $k = 1,\ldots, 10$, and $i = 1,
\ldots, \Nqoi$ we solve the optimization problem~\eqref{eq:sparse_lin_reg}; in
our computations, we use the solver SPGL1~\cite{PCE07}.  For the components of
that data vector of $\vec d$ in~\eqref{eq:sparse_lin_reg}, we use the training
set of $F_i^k$.  Therefore, every combination of $k$, $\Nord$ and $\tau$
results in a surrogate model denoted as $g^k_{\Nord, \tau}(s,\vec{\xi}^r) $.

To assess the accuracy of each bispectral surrogate we compute the average relative error
\begin{equation}\label{eq:erel}
  e_\text{rel}(g^k_{\Nord, \tau}) \!=\! 
\displaystyle\left[\frac{ \sum_{j=1}^M \int_\X \big[S(t, \vec\xi_j) - g^k_{\Nord, \tau}(s,\vec{\xi}_j^r)\big]^2
    \, ds}
{\sum_{j=1}^M \int_\X S(t, \vec\xi_j)^2 \, ds}\right]^{\Scale[1.125]{\frac12}},
\end{equation}
where $\X = [0, T_f]$, $M = 35$ and $\vec\xi_j$ is the input parameter in the
full space corresponding to $\vec{\xi}_j^r$ in the validation set of $W^k$.

We repeat the process for each of the 10 partitions, and compute the average of
$\erel$ across all partitions
\begin{equation*}\label{eq:nord_err}
  e_{\tau}^{\Nord} = \frac{1}{10} \sum_{k=1}^{10} e_\text{rel}(g^k_{\Nord, \tau}).
\end{equation*}
The cross-validation errors corresponding to $S(t, \vec\xi)$ are
displayed in Figure~\ref{fig:cross_validation_gas_saturation}. The smallest
$e_{\tau}^{\Nord}$ corresponds with $\Nord =2$ and $\tau =3.5$.
%
%------------------------------- fig ----------------------------------------
\begin{figure*}[htp!]\centering
\includegraphics[width=\textwidth]{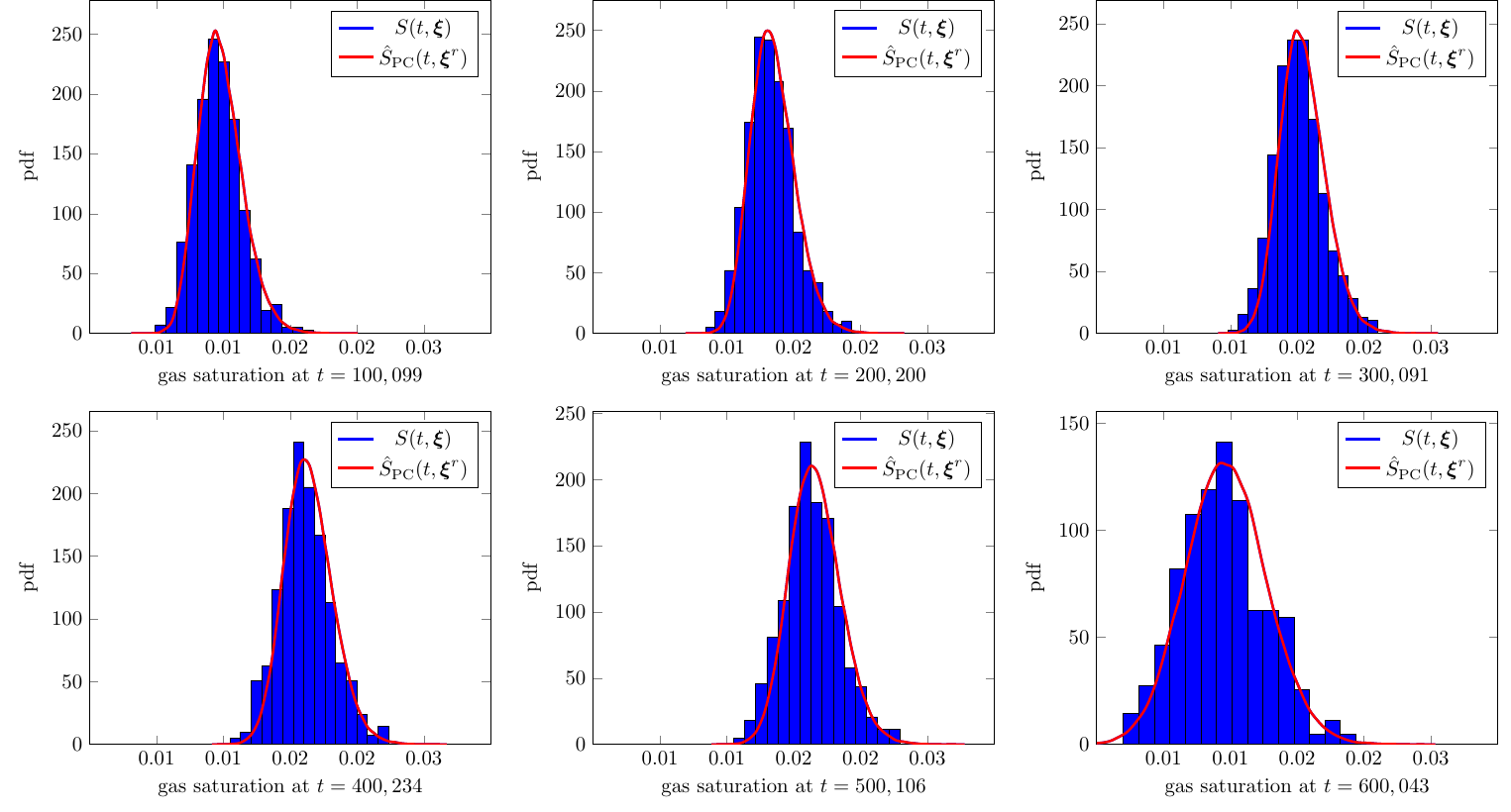}
\caption{Comparison of normalized histograms for $S(t, \vec\xi)$
 and pdf estimates of the surrogate $\SPC(t, \vec\xi^r)$ for a variety of times
$t\in [0, T_f]$.}
\label{fig:kde_gsat}
\end{figure*}

\textbf{Computing the overall bispectral surrogate.}
Once we have determined appropriate values for $\Nord$ and $\tau$ we follow
Algorithm~\ref{alg:KLE_surrogate} to construct a surrogate model from the
truncated KLE expansion of the function-valued QoI. To determine PCE 
for each KL mode $S_i(\vec\xi)$, $i=1, \ldots, \Nqoi$, we use the
solver SPGL~\cite{PCE07} to implement sparse linear regression over the \emph{entire}
350 point data set $F_i$. We use the resulting expansions to form the overall
bispectral surrogate:
\begin{equation*}\label{eq:KLE_PCE}
  \SPC = \bar{S}(t) + \sum_{i=1}^{\Nqoi} \sqrt{\lambda(\Cqoi)}
S_i^{\PC}(\vec\xi^r)   \Phi_i(t).
\end{equation*}
Note that in numerical computations, $\bar{S}(t)$ is the sample mean
$\bar{S}(t) = \frac1\Ns\sum_{j=1}^{\Ns} S(t, \vec\xi_j)$.

Next, we assess the  effectiveness of the bispectral surrogate to reflect the statistical
properties of the true model. First, we compare the sample standard deviations of
$\SPC(t, \vec\xi^r)$ and $S(t, \vec\xi)$ computed over the testing set
$\hat{W}$. The results are shown in
Figure~\ref{fig:cross_validation_gas_saturation}~(bottom). Note, the surrogate
model does an excellent job capturing the behavior of $S(t, \vec\xi)$.
Then, we compute the pdf of $\SPC(t, \vec\xi^r)$ with $100{,}000$ surrogate
evaluations and compare with the normalized histograms of the $550$ exact model
evaluations.  In Figure~\ref{fig:kde_gsat} clockwise from upper left we show
the pdf estimates for a few representative simulation times.  Note that pdf
estimates closely match the distribution of the full model.

%
%                              SECTION
%
\subsection{Gas flux at the outflow boundary}\label{sec:gas_flux_surr}
\begin{figure}[h!]
\includegraphics[width=0.45\textwidth]{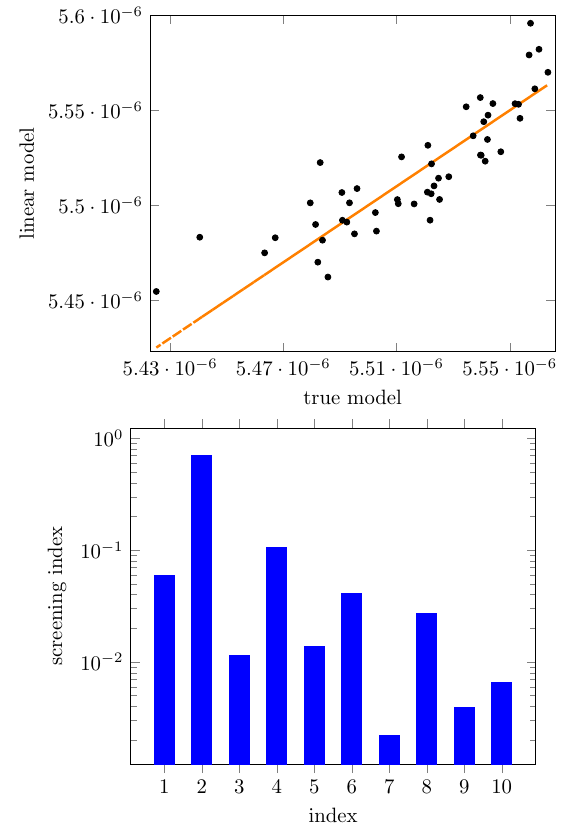}
\caption{Top: comparison of the true model to the linear model for $Q(t,\vec\xi)$
at $t=500{,}106$ years, bottom: screening indices for $Q(t, \vec\xi)$.}
\label{fig:screening_gas_flux}
\end{figure}

In this section, we study 
gas flux at the outflow boundary, denoted by $Q(t, \vec{\xi})$.
A few realizations of $Q(t, \vec{\xi})$ are shown in
Figure~\ref{fig:QoIs}~(bottom).
The global linear model is computed with $500$ model realizations.
A representation of the linear model at time $t=500{,}106$ years is displayed in
Figure~\ref{fig:screening_gas_flux}~(top). Next, we compute the screening indices
$\sj$. In Figure~\ref{fig:screening_gas_flux}~(bottom) we display
$\sj$, $j = 1\ldots, 10$ above $tol = 0.02$ only. Therefore, dimension
reduction results in the reduced input parameter
$\vec{\xi}^r = [\begin{matrix} \xi_1 & \ldots & \xi_{10}\end{matrix}]^\top$. 
Next, we compute the KLE and truncate at
$\Nqoi = 7$ terms. Then, we construct the surrogate model
using the data sets $W$ and $F_i$, where the $F_i$'s for this instance consist
of the KL modes computed for $Q(t, \vec\xi)$. We use the 10-fold cross-validation
technique described in Section~\ref{sec:gas_saturation_surr} to choose the sparse
linear regression parameters $\Nord=2$ and $\tau=4$. Finally, we use these
values to generate the bispectral surrogate $\QPC(t, \vec{\xi}^r)$.
%
%------------------------------- fig ----------------------------------------
\begin{figure*}[htp!]\centering
\includegraphics[width=\textwidth]{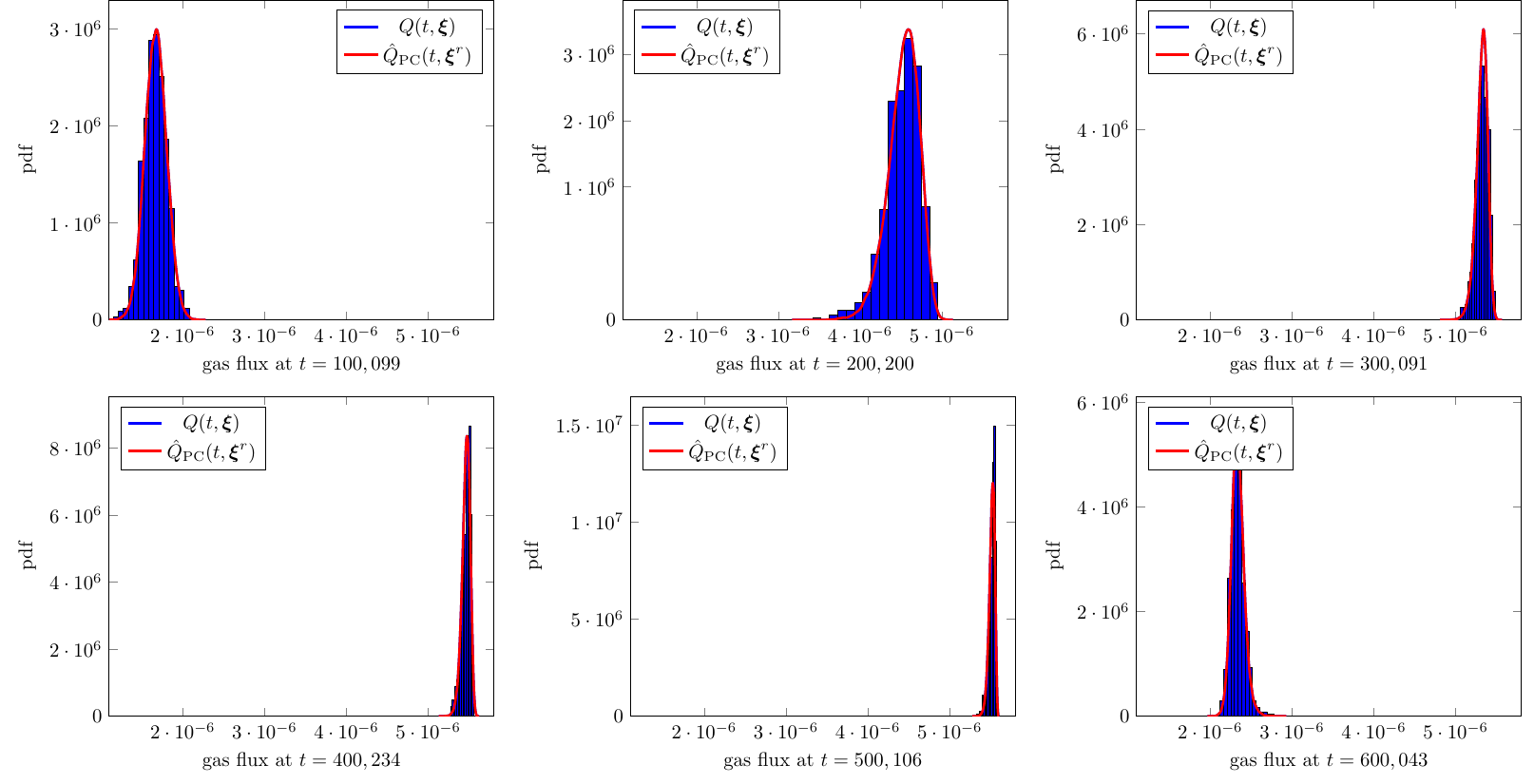}
\caption{Comparison of normalized histograms for $Q(t, \vec\xi)$
 and pdf estimates of the surrogate $\QPC(t, \vec\xi^r)$ for a variety of times
$t\in [0, T_f]$.}
\label{fig:kde_gflux}
\end{figure*}

\begin{figure}[h!]
\includegraphics[width=0.45\textwidth]{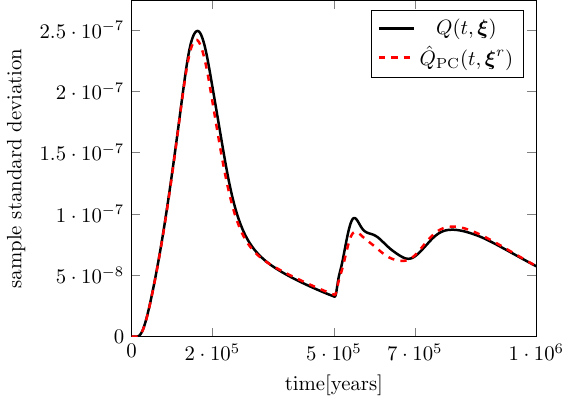}
\caption{Comparison of sample standard deviations of $Q(t, \vec\xi)$ and
$\QPC(t, \vec\xi^r)$ computed with $200$ sample points.}
\label{fig:std_dev_200_gas_flux}
\end{figure}

As before, to assess the effectiveness of the surrogate to capture the statistical
properties of the true model we compare the sample standard
deviation of the full model $Q(t, \vec{\xi})$ and the surrogate
$\QPC(t, \vec\xi^r)$, computed on $200$ validation samples. Results are displayed in
Figure~\ref{fig:std_dev_200_gas_flux}~(bottom). Lastly, using $100{,}000$
samples of $\QPC(t, \vec\xi^r)$  we compute pdf estimates at equally spaced
points in time and compare to normalized histograms created with $550$ full model
evaluations; see Figure~\ref{fig:kde_gflux}. The results in
Figure~\ref{fig:kde_gflux} and Figure~\ref{fig:std_dev_200_gas_flux} demonstrate that the
constructed surrogate for gas flux approximates the distribution of the full model reliably.

%------------------------------- fig ----------------------------------------
%
%                           SECTION
%
\subsection{Gas saturation across the domain}\label{sec:spatial_surr}

\begin{figure*}[htp!]\centering
\includegraphics[width=\textwidth]{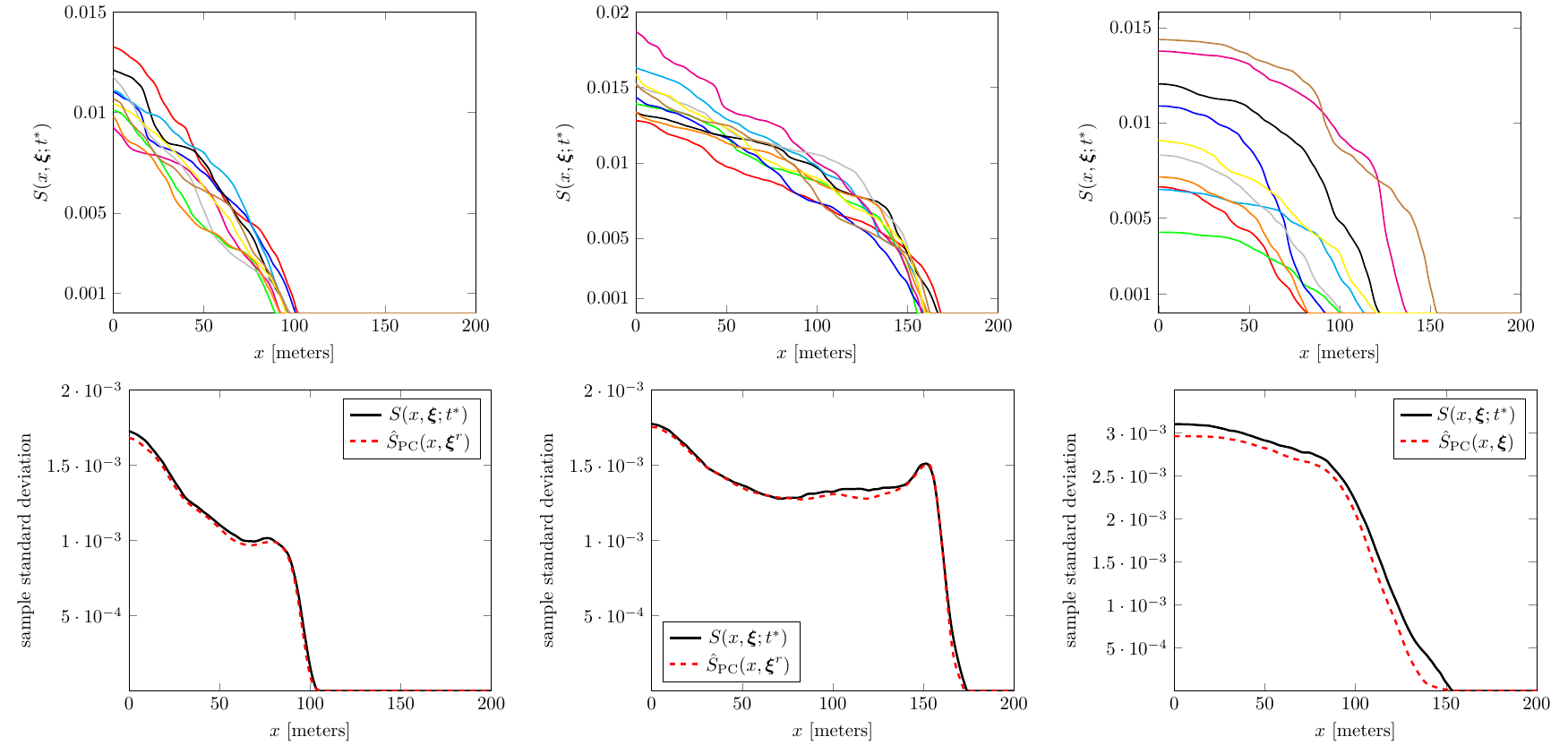}
\caption{Top row, left to right: sample realizations of $S(x, \vec{\xi}; t^*)$ for
times $100{,}099$, $300{,}091$, and $600{,}043$ years; 
bottom row, left to right:
comparison of sample standard deviation of $S(x, \vec\xi; t^*)$ and $\SPC(x, \vec\xi; t^*)$
computed on $200$ sample points.}
\label{fig:std_dev_spatial}
\end{figure*}

In this section, we focus on 
a \emph{spatially} varying QoI. Let $S(x,
\vec{\xi}; t^*)$ represent the QoI  gas saturation across the spatial domain for a
fixed time $t^*$. In particular, we include surrogate results at 
$t^* \in \{100{,}099 \bigcomma \, 300{,}091\bigcomma \, 600{,}043\}$ years.
We display several realizations for each QoI in
Figure~\ref{fig:std_dev_spatial}~(top).
The surrogate models for spatial QoIs are computed
via a similar procedure. Hence, for brevity, we include procedure
details for $t^*=600{,}043$ years only.  The relevant parameter values for the
other QoIs are included in Table~\ref{tab:surr_params}.

\begin{figure}[h!]
%------------------------------- fig ----------------------------------------
\includegraphics[width=0.45\textwidth]{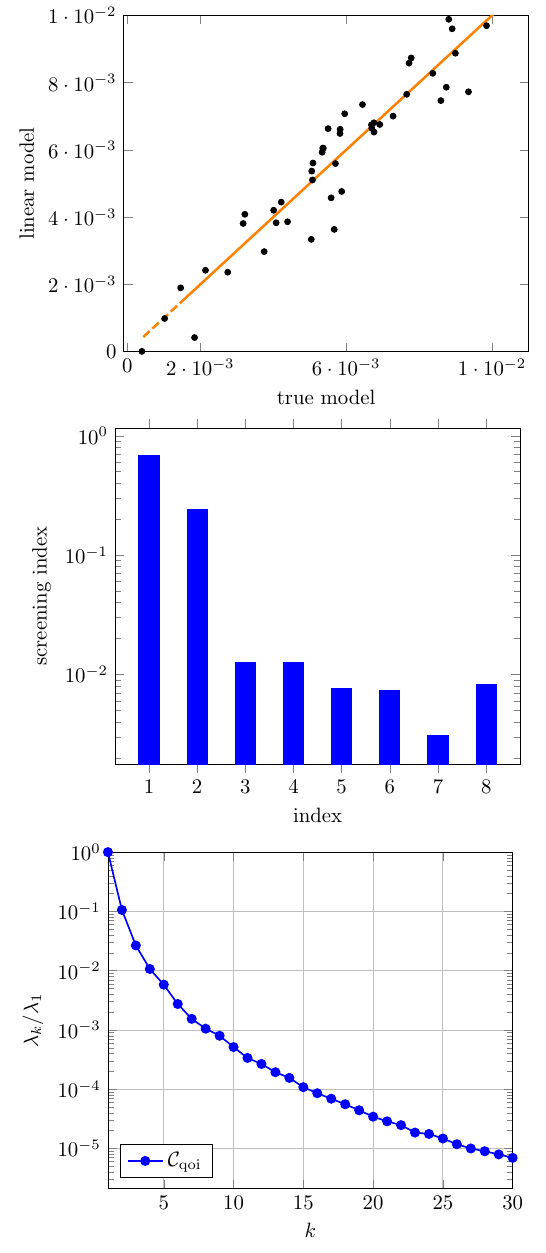}
\caption{
Results for $S(x, \vec\xi; t^*)$ with $t^*=600{,}043$ years. 
Top: comparison of the true model and the linear model for
gas saturation across the domain, middle: screening indices
for $S(x, \vec{\xi}; t^*)$, bottom:  ratio $\lambda_k/\lambda_1$, $k=1, \ldots, 30$ for
$\lambda_k(\Cqoi)$ corresponding to $S(x, \vec{\xi}; t^*)$.}
\label{fig:spatial_screening}
\end{figure}

We consider the (spatial) global linear model 
for $S(x, \vec\xi; t^*)$. 
In Figure~\ref{fig:spatial_screening}~(top) the linear model at
$x=65.5$~meters is displayed. The global linear model was observed to 
perform similarly at other
values of $x$.  Next, we compute the screening indices $\sj$ and use the
importance tolerance $tol = 0.002$ for dimension reduction resulting in the
reduced input parameter $\vec{\xi}^r = [
\begin{matrix} \xi_1 & \xi_2 & \ldots & \xi_8\end{matrix}]^\top$.  In
Figure~\ref{fig:spatial_screening}~(middle) we display the screening indices
corresponding to these parameters.

Next, we compute the KLE of $S(x, \vec\xi; t^*)$ 
using Nystr\"{o}m's method
with $550$ model evaluations. In Figure~\ref{fig:spatial_screening}~(bottom) we
report the normalized eigenvalues of the output covariance operator $\Cqoi$ for
$S(x, \vec\xi; t^*)$. This result is included to demonstrate that the gas
saturation process is also low-rank in space. We truncate the KLE at $\Nqoi =
5$ terms.

As before, the PCE for the KL modes are computed with
sparse linear regression using 350 full model realizations. Once again, the
cross-validation procedure described in Section~\ref{sec:gas_saturation_surr} is used
to determine $\Nord = 3$ and $\tau=2.8$. Lastly, the computed PCEs for each KL mode
is used to construct the bispectral surrogate $\SPC(x, \vec\xi^r)$.

\begin{table}[h!]
\centering
\renewcommand{\arraystretch}{1.5}
\begin{tabular}{|c|c|c|c|c|}
\hline
surrogate          & for fixed $t$ or $x$ & $\Nqoi$ & $\Nord$ &       error       \\
\hline
$\SPC(t, \vec\xi)$ & $x = 0$ meters       & $5$     & 2       & $3.4813\cdot10^{-2}$\\
\hline
$\QPC(t, \vec\xi)$ & $x = 200$ meters     & $7$     & 2       & $7.5019\cdot10^{-3}$ \\
\hline
$\SPC(x, \vec\xi)$ & $t^*=100{,}099$ years & $7$     & 2       & $3.0397\cdot10^{-2}$ \\
\hline
$\SPC(x, \vec\xi)$ & $t^*=300{,}091$ years & $11$    & 2       & $2.1690\cdot10^{-2}$ \\
\hline
$\SPC(x, \vec\xi)$ & $t^*=600{,}043$ years  & $5$     & 3       & $8.3110\cdot10^{-2}$ \\
\hline
\end{tabular}
\caption{Surrogate parameter values  and $\erel$ errors for surrogate models.}
\label{tab:surr_params}
\end{table}

To evaluate the effectiveness of the surrogate models for 
$t^* \in \{ 100{,}099\bigcomma \, 300{,}091\bigcomma \, 600{,}043\}$, 
we compare the sample
standard deviation of $S(x, \vec\xi; t^*)$ and $\SPC(x, \vec\xi^r)$ for 200 sample
points. These results are displayed in
Figure~\ref{fig:std_dev_spatial}~(bottom).
Observe, for $t^*=100{,}099$ and $t^*=300{,}091$ years the surrogate model
replicates the sample standard deviation well. For $t^*=600{,}043$ years note that
while we are underestimating the sample standard deviation, we are
still capture the overall behavior of the full model.
%
%
%------------------------------- fig ----------------------------------------

The capability of the computed bispectral surrogate to replicate true
model behavior can also be tested by computing  the average relative error
defined in~\eqref{eq:erel}.Table~\ref{tab:surr_params} contains the values
for $\erel$ computed over the validation set $\hat{W}$ for each surrogate
presented in this section, as well as those in
Sections~\ref{sec:gas_saturation_surr} and \ref{sec:gas_flux_surr}. Note
that for the spatially varying QoIs, we let $\X = [0,200]$ and for temporally
varying we let $\X = [0, T_f]$, in~\eqref{eq:erel}. Note, the error across all
surrogates is less than 8\%, and in four out the five surrogates is less than
$4$\%.  The largest $\erel$ corresponds to $\SPC(c, \vec\xi)$ at $t=600{,}043$
years, in which case we are also underestimating the standard deviation.

%
%                             SECTION
%
\subsection{Using the surrogate model}\label{sec:use_surrogate}
Here we illustrate the use of surrogates for temporally varying QoIs in
performing statistical studies.  In particular, we perform model prediction,
variance-based global sensitivity analysis (GSA) by computing Sobol' indices, 
and a study of output correlation
structure.  It is worth noting that Sobol' indices can be computed analytically
when using bispectral surrogates; see \cite{AlexanderianGremaudSmith20}.
However, to keep the discussion general and since the cost of evaluating the
bispectral surrogate is negligible, here we rely on sampling the bispectral
surrogate for performing GSA.  Specifically, we also perform GSA on QoIs that
are derived from the bispectral surrogates such as maximum gas saturation and
maximum gas flux at the outflow boundary.

\begin{figure}
\includegraphics[width=0.45\textwidth]{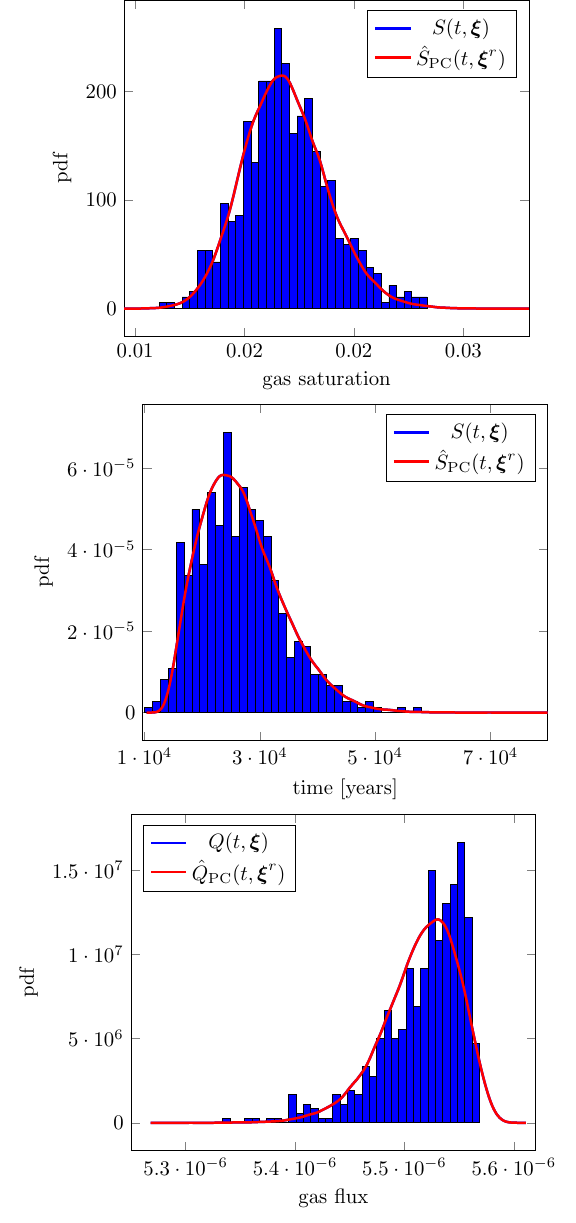}
\caption{Comparison of normalized histograms and pdf estimates for top: max saturation
value of $S_{\max}$, middle: first time $\SPC(t, \xi^r)$ is above $20$\% $S_{\max}$,
bottom: max flux value $Q_{\max}$.}
\label{fig:surrogate_prediction}
\end{figure}

\textbf{Model prediction.}
We consider using $\SPC(t, \vec\xi^r)$ and \\
$\QPC(t, \vec\xi^r)$ for making
predictions.  Recall, these bispectral surrogates correspond to gas saturation at the
inflow boundary and gas flux at the outflow boundary. We study three
observables of interest: maximum gas saturation, denoted $S_{\max}$, maximum
gas flux, denoted $Q_{\max}$, and the first time for which gas saturation rises
above 20\% of $S_{\max}$. We compute $100{,}000$ realizations of each
surrogate, extract the pertinent observables, and use the samples to compute
pdf estimates.  In Figure~\ref{fig:surrogate_prediction}, we compare the pdf
estimates against the normalized histograms computed using exact model
evaluations.  These results indicate the utility of the surrogates for
estimating the statistical properties of model observables.
%
%------------------------------- fig ----------------------------------------
%%
%------------------------------- fig ----------------------------------------
%
\begin{figure}
\includegraphics[width=0.45\textwidth]{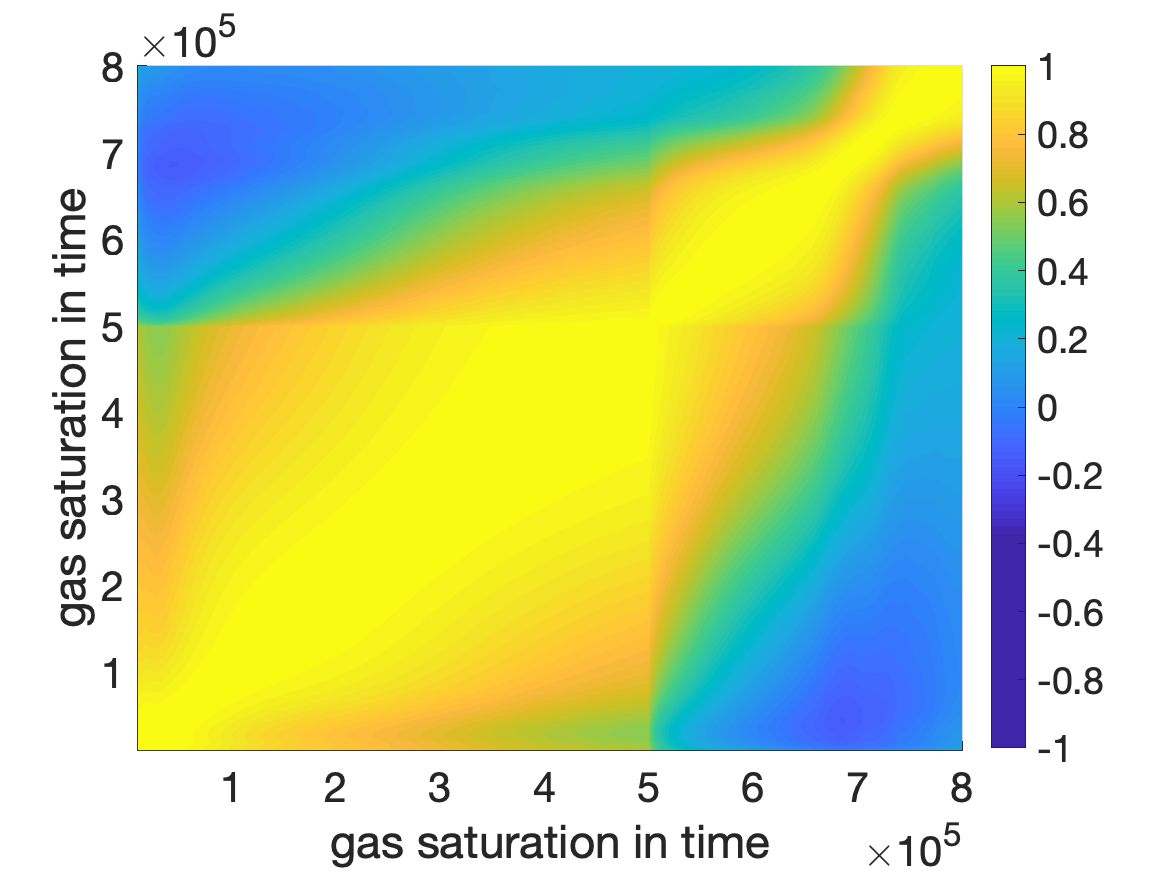}
%------------------------------- fig ----------------------------------------
\includegraphics[width=0.45\textwidth]{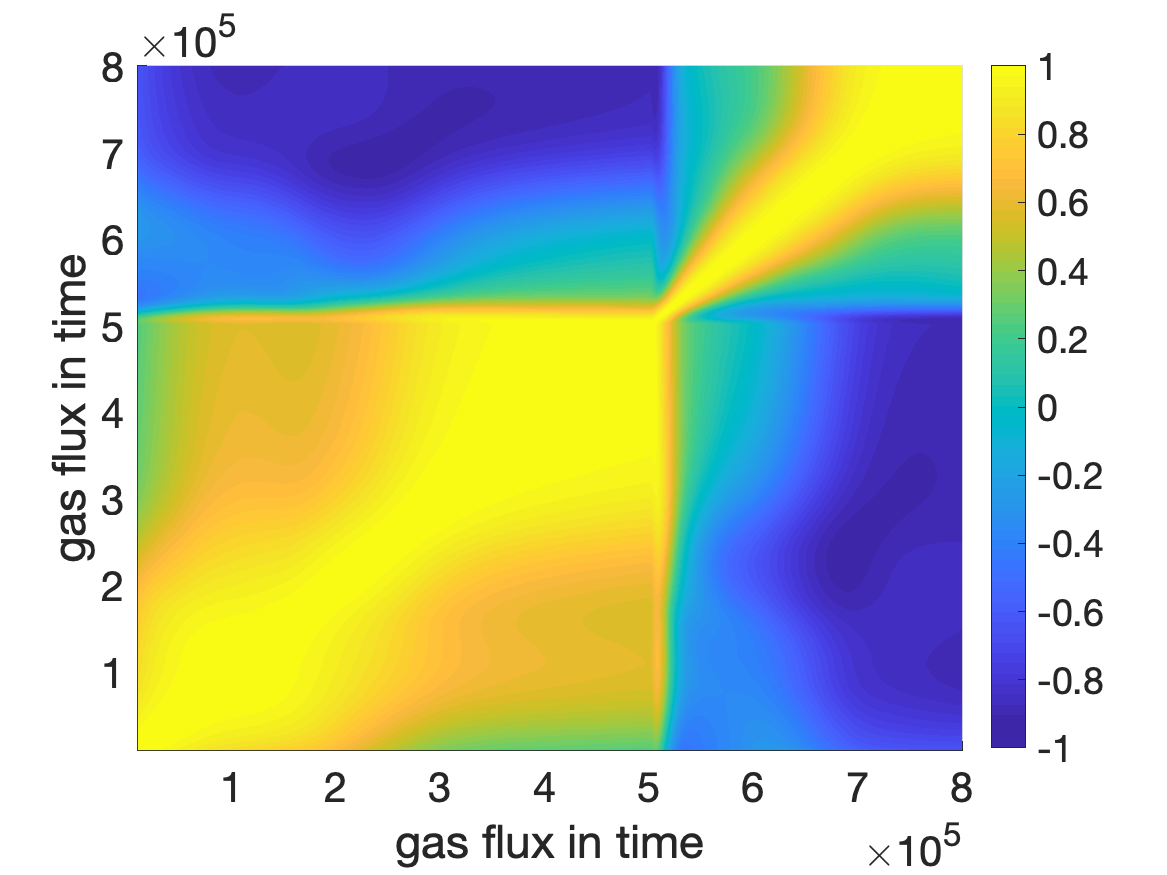}
%------------------------------- fig ----------------------------------------
\includegraphics[width=0.45\textwidth]{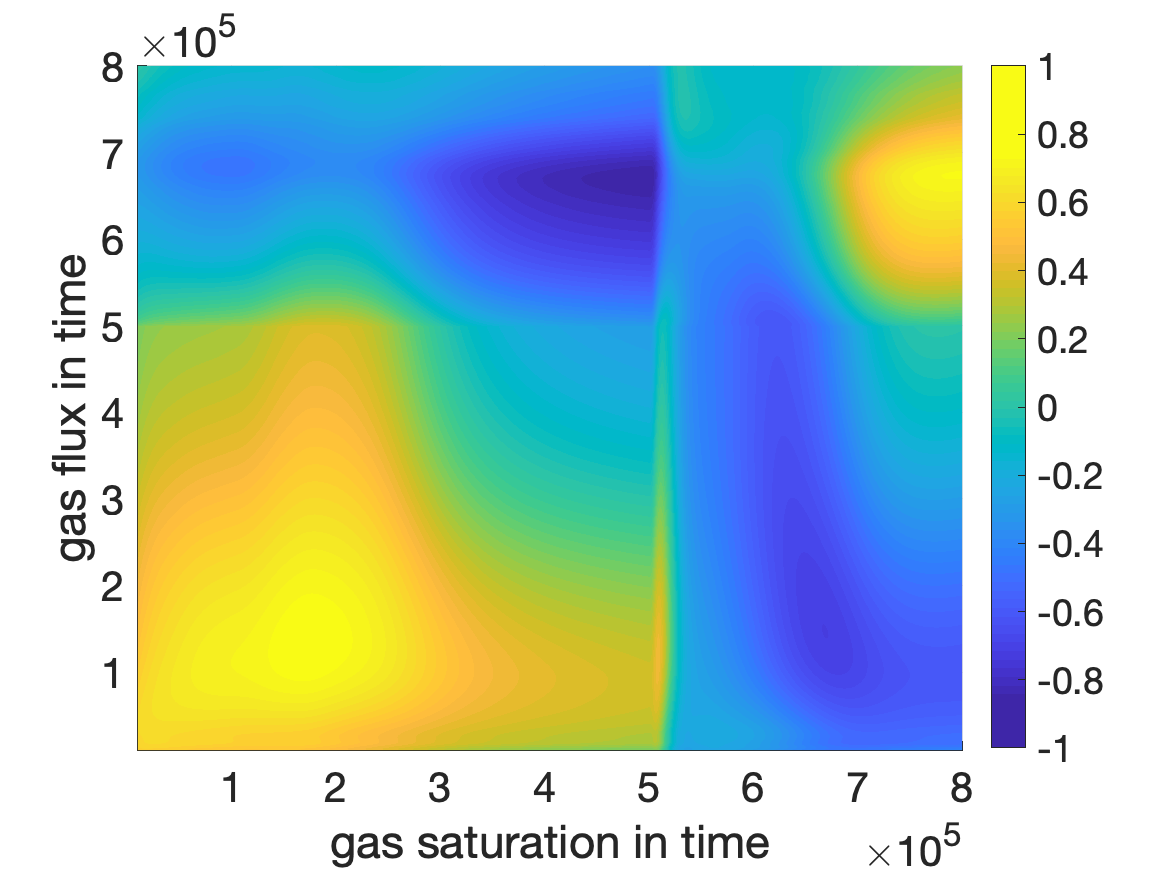}
%------------------------------- fig ----------------------------------------
\caption{Top: correlation matrix for $\SPC(t, \vec\xi^r)$ computed using
the analytic formula in~\eqref{eq:corr}, middle: correlation matrix for
$\QPC(t, \vec\xi^r)$ computed using the analytic formula in~\eqref{eq:corr}
bottom: cross-correlation structure of $\SPC(t, \vec\xi^r)$ and
$\QPC(t, \vec\xi^r)$ computed using the analytic formula in~\eqref{eq:crosscorr}.
}
\label{fig:corr}
\end{figure}

\textbf{Variance based sensitivity analysis via Sobol' indices.}
total Sobol' indices provide an informative global sensitivity analysis
tool that apportions percentages of QoI variance due to input parameter
variations. While total Sobol' indices are traditionally applied to scalar
QoIs~\cite{Sobol01,Sobol90}, there exist extensions for variance based
analysis to function-valued
QoIs~\cite{GamboaJanonKlein14,AlexanderianGremaudSmith20}, referred to as
\emph{functional} total Sobol' indices.

In general, calculating Sobol' indices for computationally intensive models is
challenging. This involves an expensive sampling procedure that requires
a large number of model evaluations. An efficient-to-evaluate surrogate model can be
used to accelerate this process.  We use the temporal surrogates to compute
total Sobol' indices  for both function-valued and scalar QoIs. In particular,
we compute the functional total Sobol' indices for $\SPC(t,\vec\xi^r)$ and
$\QPC(t, \vec\xi^r)$, both of which are functions in $t$, and we compute the
total Sobol' indices for the scalar QoIs $S_{\max}$ and $Q_{\max}$.  In each
case, we compute the total Sobol' indices via sampling, using a variety of
samples sizes: $\Ns=\{1{,}000\bigcomma\, 10{,}000\bigcomma\, 50{,}000\}$.

\begin{figure*}
\begin{tabular}{cc}
\includegraphics[width=0.45\textwidth]{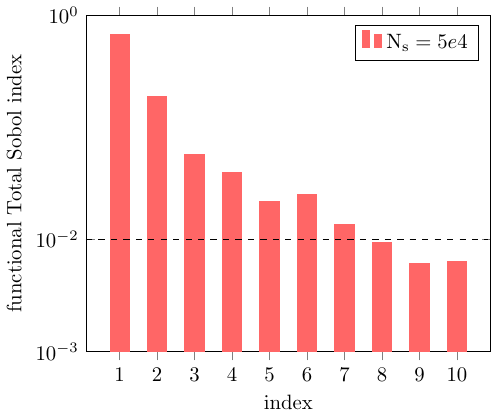} &
% ------------------------------- fig ----------------------------------------
\includegraphics[width=0.45\textwidth]{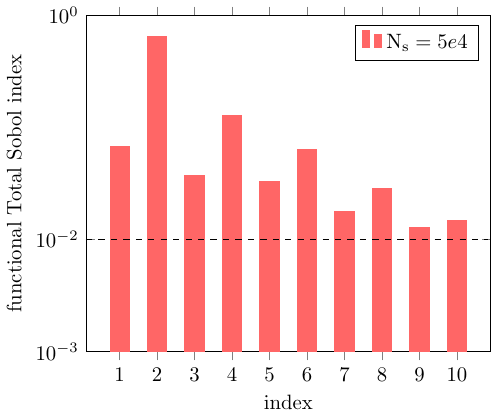} \\
% ------------------------------- fig ----------------------------------------
\includegraphics[width=0.45\textwidth]{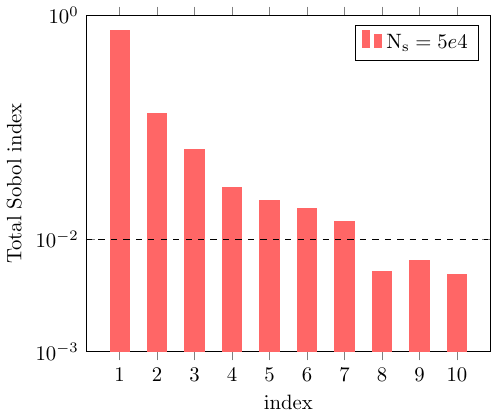}&
% ------------------------------- fig ----------------------------------------
\includegraphics[width=0.45\textwidth]{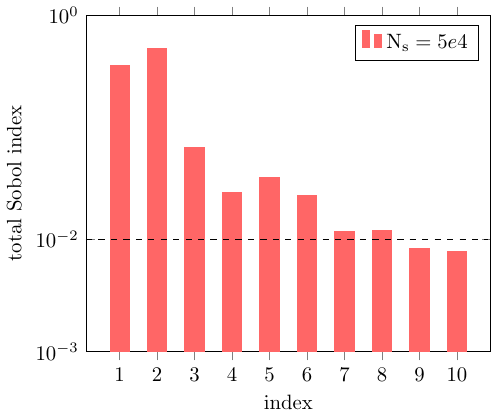}\\
\end{tabular}
\caption{From top left counter clockwise: functional total Sobol' indices for
$\SPC(t, \vec\xi^r)$, functional total Sobol' indices for
$\QPC(t, \vec\xi^r)$, total Sobol' indices for $S_{\max}$, total Sobol' indices for
$Q_{\max}$.}
\label{fig:Tsobol}
\end{figure*}

The results in the top row of Figure~\ref{fig:Tsobol} show the functional
Sobol' indices for $\SPC(t,\vec\xi^r)$ and $\QPC(t, \vec\xi^r)$.
Note that the magnitudes in the top row of Figure~\ref{fig:Tsobol} are similar
to those in Figure~\ref{fig:screening_gas_saturation} and
Figure~\ref{fig:screening_gas_flux}.  This provides further support for the
original input parameter importance ranking and subsequent dimension reduction.
In the bottom row of Figure~\ref{fig:Tsobol} we report the total Sobol' indices
for $S_{\max}$ and $Q_{\max}$. We also note that for the gas saturation
QoIs~(Figure~\ref{fig:Tsobol}(left: top and bottom), the importance ranking of
the input parameters is similar. In contrast, there is more variability in
ranking for gas flux QoIs~(Figure~\ref{fig:Tsobol} (right: top and bottom).

Finally, we mention that for many applications, the total Sobol' indices can be
used for further input parameter dimension reduction. For the present model
however, we did not reduce the input parameter further because the surrogate
model computed was already efficient and sufficiently accurate.

%
% ------------------------------- fig ----------------------------------------

\textbf{Correlation structure}
Lastly, we illustrate the use of the bispectral surrogates for computing the correlation
structure of the output, which is a useful tool for understanding overall model
dynamics. Using equation~\eqref{eq:corr} we compute the correlation function of
$\SPC(t, \vec{\xi}^r)$ and $\QPC(t, \vec\xi^r)$.  The resulting heat maps are
shown in Figure~\ref{fig:corr} top and middle, respectively. The results for
$\SPC(t, \vec\xi^r)$ suggest significant correlations across time. This
behavior is also seen in  the correlation function of $\QPC(t, \vec\xi^r)$,
except the sudden shift in dynamics at the time $t=500{,}000$ years; recall,
this the time gas injection stops.
We also compute the cross-correlation between $\SPC(t, \vec\xi^r)$ and
$\QPC(t, \vec\xi^r)$ using the formula in~\eqref{eq:crosscorr};
see Figure~\ref{fig:corr}~(bottom).
The heat map  suggests there is large cross-correlation  between the two QoI
for both early and late times.

\section{Conclusion}
\label{sec:conc}
We have presented a structure exploiting non-intrusive framework for efficient dimension
reduction and surrogate modeling for models with high-dimensional
inputs and outputs. The proposed parameter screening metric utilizes
approximate global sensitivity measures for function-valued outputs that rely
on concepts from global sensitivity analysis and active subspace methods.
An efficient bispectral surrogate model was constructed from a truncated KLE of the QoI by
approximating the KL modes with PCEs. Note, these KL mode PCEs were constructed
in the reduced parameter space. 

We deployed our framework for fast uncertainty analysis in a multiphase
multicomponent flow model. The efficiency and effectiveness of the surrogate
model was demonstrated  with a comprehensive set of numerical experiments, where
we consider a number of function-valued (temporally or spatially distributed)
QoIs.  In particular, our results indicate that it is possible to use a modest
amount of model realizations to reduce both the input and output dimensions and
construct an efficient surrogate model. The proposed framework not only
provides efficient surrogates, it also reveals and exploits the
low-dimensional structures in model input and output spaces, which provides
further insight into the behavior of the governing model.

In general, the screening approach takes the following form. We construct a
cheap approximation $\tilde{f}$, compute the corresponding screening
indices~\eqref{eq:screening_indc}, and use them to reduce the input parameter
space. Our approach relies on the screening metrics being sufficiently accurate
surrogates and cheap to compute for the derivative-based global sensitivity
measures for the function-valued QoIs under study.  This in turn assumes the
global linear model constructed within the parameter screening procedure leads
to a sufficient approximation of the activity scores.  It is observed that this
global linear model can successfully capture one-dimensional active subspaces
in a wide range of applications~\cite{Constantine15}. The success of this
strategy for obtaining approximate activity scores was also observed in the
present work, in the context of a complex nonlinear flow model. However, for
models that exhibit highly nonlinear parameter dependence a linear model might
fail to provide accurate global sensitivity information.
In~\cite{GreyConstantine18,SeshadriSeshadriConstantineEtAl18}, global quadratic
models were used effectively to accelerate active subspace discovery for
scalar-valued QoIs. Exploring quadratic models within our framework provides an
interesting direction for future work and would allow application of the
proposed strategy to a broader class of problems.

\begin{acknowledgements}
The research of  H. Cleaves and A. Alexanderian was partially supported
by the National Science Foundation through the grant DMS-$1745654$.
The work of A.~Alexanderian was also supported in part through 
the grant DMS-$1953271$.
\end{acknowledgements}

% Authors must disclose all relationships or interests that
% could have direct or potential influence or impart bias on
% the work:
%
%\section*{Conflict of interest}
%The authors declare that they have no conflict of interest.
\section*{Data availability}
Data sharing is not applicable to this article as no datasets were generated or
analysed during the current study.
The results reported in this article are all based on numerical simulations.
The details of the benchmark problem considered can be found in~\cite{Bourgeat2012}.

\section{Appendix}
\label{sec:append}
\subsection{Proof of upper bound on total error in product space}
\begin{proof}
Let $f(s, \vec\xi)$ be in $L^2$ of the product space $\paramspace\times\X$ and
$\norm{\cdot}$ be the $L^2$ error in the product space $\paramspace\times\X$. The
truncated KLE of $f$ is given by
\[
  \fPC(s, \vec\xi) = \bar{f}(s) +
  \sum_{i=1}^{\Nqoi} \sqrt{\lambda_i} f_i^{\PC}(\vec\xi) \Phi_i(s).
\]
The total error in the product space is given by
\[
  \norm{f-\fPC}^2 \leq 2\norm{f-\fn}^2 + 2\norm{\fn - \fPC}^2
\]
We consider the first term
\begin{align*}
  & 2\norm{f-\fn}^2\\
  =& ~2\norm{\sum_{i=1}^{\infty}\sqrt{\lambda_i}f_i(\vec\xi)\Phi_i(s) - 
  \sum_{i=1}^{\Nqoi}\sqrt{\lambda_i}f_i(\vec\xi)\Phi_i(s)}^2\\
  =& ~2\hspace{-1mm}\int_\paramspace\int_\X \left(\sum_{i=\Nqoi+1}^{\infty}\sqrt{\lambda_i}
  f_i(\vec\xi)\Phi_i(s)\right)^2 \,ds\,\mymu(d\vec\xi)\\
  =& ~2\hspace{-3mm}\sum_{i,j =\Nqoi+1}^\infty \hspace{-2mm}\sqrt{\lambda_i}\sqrt{\lambda_j}  \int_\paramspace
  f_i(\vec\xi)f_j(\vec\xi)\int_\X \Phi_i(s)\Phi_j(s)\,ds \,\mymu (d\vec\xi)\\
  =& ~2\hspace{-3mm}\sum_{i=\Nqoi+1}^\infty\hspace{-2mm}\lambda_i \int_\paramspace f_i(\vec\xi)^2\,\mymu (d\vec\xi)
  = ~2\hspace{-3mm}\sum_{i=\Nqoi+1}^\infty\hspace{-2mm} \lambda_i.
\end{align*}
Changing the order of infinite sums and integral is a consequence of the Dominated
Convergence Theorem and reordering of integrals is a justified by Fubini's
Theorem. The orthogonality of the eigenfunctions in $L^2(\X)$ justifies the
simplification in the second to last line, and the last step is a consequence of
the KL modes properties.

Next, we consider the second error term.
Let 
\[
f_i^{\PC} = \sum_{k=0}^{\Npc} \hat{c}_{i,k} \Psi_k(\vec\xi),
\] 
where
$\hat{c}_{i,k}$ represents the numerical approximation of the exact PCE
coefficients $c_{i,k}$ and recall,
$f_i = \sum_{k=0}^{\infty} c_{i,k}\Psi_k(\vec\xi)$ we have
\begin{align*}
  &~ 2\norm{\fn - \fPC}^2 = 2\norm{\sum_{i = 1}^{\Nqoi} \sqrt{\lambda_i}f_i(\vec\xi)\Phi_i(s) 
  -\sum_{i=1}^{\Nqoi}\sqrt{\lambda_i}f_i^{\PC}(\vec\xi) \Phi(s)}^2\\
  =&~ 2\int_\paramspace \int_\X \left(\sum_{i=1}^{\Nqoi} \sqrt{\lambda_i}\Phi_i(s)\left[f_i(\vec\xi)-f_i^{\PC}(\vec\xi)\right]\right)^2\,ds\,\mymu (d\vec\xi)\\
  =&~ 2\sum_{i,j = 1}^{\Nqoi} \sqrt{\lambda_i}\sqrt{\lambda_j} \int_\paramspace
  (f_i-f_i^{\PC})(f_j-f_j^{\PC}) \int_\X \Phi_i(s)\Phi_j(s)\,ds\,\mymu (d\vec\xi)\\
  =&~ 2\sum_{i=1}^{\Nqoi} \lambda_i\int_\paramspace (f_i(\vec\xi)-f_i^{\PC}(\vec\xi))^2\,
  \mymu(d\vec\xi)\\
  =&~ 2\sum_{i=1}^{\Nqoi}\lambda_i \int_\paramspace \left(\sum_{k =0}^\infty c_{i,k}\Psi_k(\vec\xi) - \sum_{k=0}^{\Npc}\hat{c}_{i,k}\Psi_k(\vec\xi)\right)^2\,\,\mymu (d\vec\xi)\\
  =&~ 2\sum_{i = 1}^{\Nqoi} \lambda_i\int_\paramspace \left(\sum_{k=0}^{\Npc} (c_{i,k}-\hat{c}_{i,k})\Psi_k(\vec\xi) +\hspace{-2.0mm}\sum_{k=1+\Npc}^\infty \hspace{-3.0mm}c_{i,k}\Psi_k(\vec\xi)\right)^2\,\mymu (d\vec\xi)\\
  =&~ 2\sum_{i=1}^{\Nqoi} \lambda_i \sum_{k=1}^{\Npc}(c_{i,k}-\hat{c}_{i,k})^2
  \normLomega{\Psi_k}^2\\
  &\hspace{1.5cm}+~  2\sum_{i=1}^{\Nqoi} \lambda_i \sum_{j=1+\Npc}^{\infty}
  c_{i,j}^2\normLomega{\Psi_j}^2.
\end{align*}
The simplification in the third line a consequence of the orthogonality of the
PCE basis functions.

Thus, we have a bound on the total error
\begin{align*}
  & \norm{f-\fPC}^2 \leq 2\norm{f-\fn}^2 + 2\norm{\fn - \fPC}^2\\
  &=~ 2\sum_{i=\Nqoi+1}^\infty\hspace{-2mm}\lambda_i + \sum_{i=1}^{\Nqoi} \lambda_i  \sum_{k=1}^{\Npc}(c_{i,k}-\hat{c}_{i,k})^2
\normLomega{\Psi_k}^2\\
&+  2\sum_{i=1}^{\Nqoi} \lambda_i \sum_{j=1+\Npc}^{\infty}
  c_{i,j}^2\normLomega{\Psi_j}^2.~\square
\end{align*}
\end{proof}

% BibTeX users please use one of
%\bibliographystyle{spbasic}      % basic style, author-year citations
%\bibliographystyle{spmpsci}      % mathematics and physical sciences
%\bibliographystyle{spphys}       % APS-like style for physics
%\bibliography{}   % name your BibTeX data base

% Non-BibTeX users please use
\bibliographystyle{plain}
\bibliography{biblio}
%
% and use \bibitem to create references. Consult the Instructions
% for authors for reference list style.

\end{document}